\documentclass[11pt]{article}
\usepackage{amsmath}
\usepackage{amsfonts}
\usepackage[a4paper]{geometry}
\usepackage{graphicx}
\usepackage[onehalfspacing]{setspace}
\usepackage{caption}
\usepackage{epstopdf}

\newtheorem{theorem}{Theorem}

\newtheorem{corollary}[theorem]{Corollary}

\newenvironment{proof}[1][Proof]{\noindent\textbf{#1.} }{\ \rule{0.5em}{0.5em}}
\input{tcilatex}
\begin{document}

\title{Generalised vector products and metrical trigonometry of a tetrahedron%
}
\author{G A Notowidigdo* \\
School of Mathematics and Statistics\\
UNSW Sydney\\
Sydney, NSW, Australia\\
gnotowidigdo@zohomail.com.au \and N J Wildberger \\
School of Mathematics and Statistics\\
UNSW Sydney\\
Sydney, NSW, Australia\\
n.wildberger@unsw.edu.au}
\date{}
\maketitle

\begin{abstract}
We study the general rational trigonometry of a tetrahedron, based on
quadrances, spreads and solid spreads, using vector products associated to
an arbitrary symmetric bilinear form over a general field, not of
characteristic two. This gives us algebraic analogs of many classical
formulas, as well as new insights and results. In particular we derive
original relations for a tri-rectangular tetrahedron.
\end{abstract}



\section{Introduction}

In this paper, which is a follow-up to \textit{Generalised vector products
applied to affine and projective rational trigonometries in three dimensions}
(henceforth referred to as \cite{NotoWild1}), we will apply the framework of
generalised scalar and $B$-vector products developed there to set up a
framework for the rational trigonometry of a general tetrahedron in
three-dimensional affine space with a general metrical structure defined
over an arbitrary field, not of characteristic two.

The metrical structure of the tetrahedron is a much studied topic, since at
least the time of Tartaglia, who was the first to give a formula for the
volume in terms of the (squared) lengths of the sides, a formula also found
by Euler and generalized by Cayley and Menger. This classical treatment of
the trigonometry of the tetrahedron reached a high point with the large
scale summary of Richardson \cite{Richardson} for the Euclidean case,
involving not just lengths and angles of faces but also dihedral angles
between faces and solid angles at vertices and transversal lengths between
opposite edges. Since locally at a vertex a tetrahedron determines a tripod,
there was also a projective aspect which connects naturally with spherical
trigonometry \cite{Todhunter}. In modern times both computational geometry
in three dimensions \cite{Golias Tsiboukis} and the finite element method in
space \cite{Brandts} have utilized tetrahedral meshes and their
measurements. And of course there is naturally interest in generalizations
to metrical formulas for more general polytopes, as in the rigidity results
of Sabitov \cite{Sabitov}, and to explorations of analogs of trigonometric
concepts and laws to simplices in higher dimensions, for example \cite%
{Eriksson}. But it is clear that our understanding of the $n$-dimensional
story is still very much in its early days.

Along with this classical mostly Euclidean orientation, since the 19th
century development of non-Euclidean geometry, the question of how
trigonometry extends to hyperbolic and spherical settings has been of keen
interest, and an important issue is the three-dimensional situation for
hyperbolic tetrahedra, where even the volume formula originally due to
Lobachevsky is much more complicated than in the Euclidean case, and much
work has been done in understanding and extending this result for example
the work of \cite{Kellerhals}, \cite{Derevnin Mednykh} and \cite{Mednykh
Pashkevich}.

With the advent of rational trigonometry (see \cite{WildDP} and \cite{WildPS}%
) the possibility emerged to reframe trigonometry in a purely algebraic
fashion; where lengths and angles become secondary to purely algebraic
concepts defined in terms of a symmetric bilinear form. A projective version
of this theory applies also to allow a recasting of hyperbolic geometry (see 
\cite{WildUHG1}, \cite{WildUHG2}, \cite{WildUHG3} and \cite{AlkWild}) which
also extends the subject to general fields, including finite fields.

In \cite{NotoWild1} we laid out a three-dimensional version of this theory
which extends classical formulas of Lagrange, Binet, Cauchy and others to
this more general setting valid for arbitrary fields (not of characteristic
two) and to general quadratic forms. In this follow up paper we intend to
apply this technology to completely reformulating the classical trigonometry
of a tetrahedron, using rational analogs also of dihedral angles and solid
spreads. We hope that this will provide a bridge also to a projective
version which describes the trigonometry of spherical and hyperbolic
tetrahedra, and in fact some key projective formulas for that project are
already contained in this present work.

Suppose that $\mathbb{V}^{3}$ is the three-dimensional vector space, over a
field $\mathbb{F}$ not of characteristic $2$, consisting of row vectors $%
v=\left( x,y,z\right) $ and that $B$ is a symmetric non-degenerate $3\times
3 $ matrix, so that $\det B\neq 0$. Then we may define a symmetric bilinear
form, or $B$\textbf{-scalar product}, by the rule that for\textbf{\ }any two
vectors $v$ and $w$ 
\begin{equation*}
v\cdot _{B}w\equiv vBw^{T}.
\end{equation*}%
This number is always an element of the field $\mathbb{F}$.

The $B$\textbf{-quadrance} of a vector $v$ is%
\begin{equation*}
Q_{B}\left( v\right) \equiv v\cdot _{B}v
\end{equation*}%
and a vector $v$ is $B$\textbf{-null} precisely when $Q_{B}\left( v\right)
=0.$

If $v$ and $w$ are non-null vectors, then the $B$\textbf{-spread} between
them is the number%
\begin{equation*}
s_{B}\left( v,w\right) \equiv 1-\frac{\left( v\cdot _{B}w\right) ^{2}}{%
Q_{B}\left( v\right) Q_{B}\left( w\right) }.
\end{equation*}%
In Euclidean geometry, the quadrance and spread are typically thought of
respectively as the squared distance and the squared sine of an angle; in
our framework, we are using quadrances and spreads to allow for extensions
of Euclidean geometry to arbitrary symmetric bilinear forms over a general
field more easily.

In \cite{NotoWild1}, we extended the definition of vector products in
Euclidean three-dimensional vector space also to this more general
three-dimensional situation. Given the usual (Euclidean) vector product $%
v\times w$, for two vectors $v$ and $w$ in $\mathbb{V}^{3}$, we define the $%
B $\textbf{-vector product} of $v$ and $w$ to be%
\begin{equation*}
v\times _{B}w\equiv \left( v\times w\right) \func{adj}B
\end{equation*}%
where $\func{adj}B=\left( \det B\right) B^{-1}$ is the adjugate of the
invertible matrix $B$.

After a short review of the properties of $B$-scalar and $B$-vector products
which were proven in \cite{NotoWild1}, we define the fundamental
trigonometric invariants in three-dimensional affine space, denoted by $%
\mathbb{A}^{3}$, where $\mathbb{V}^{3}$, as given above, is its associated
vector space. These include the $B$-quadrance and $B$-spread, as well as the 
$B$\textbf{-quadrea} of a triangle in $\mathbb{A}^{3}$ which extends the
definition of quadrea in \cite{WildDP}.

While these three quantities featured prominently in \cite{NotoWild1}, this
paper introduces four new trigonometric invariants: the $B$\textbf{-quadrume}%
, the $B$\textbf{-dihedral spread}, the $B$\textbf{-solid spread} and the $B$%
\textbf{-dual solid spread}. The $B$-quadrume, which is a quadratic version
of volume in our framework, has a close connection to the Cayley-Menger
determinant, as seen in \cite{Audet}, \cite[pp. 285-289]{Dorrie} and \cite[%
pp. 124-126]{Sommerville}. The latter three quantities, which are analogs of
angles between two planes or solid angles between three lines in our
framework, have origins in projective geometry (see \cite{WildPS} and \cite%
{Wild2Planes}) and thus highlight the power of using generalised vector
products to explain affine rational trigonometry in three dimensions.

We can then compute these quantities for a general tetrahedron in $\mathbb{A}%
^{3}$, for which we can then discover interesting algebraic relations. For a
tetrahedron $\overline{A_{0}A_{1}A_{2}A_{3}}$ with points $A_{0}$, $A_{1}$, $%
A_{2}$, $A_{3}$, we will denote, for indices $0\leq i<j<k\leq 3$:

\begin{itemize}
\item by $Q_{ij}$ the $B$-quadrance between two points $A_{i}$ and $A_{j}$;

\item by $\mathcal{A}_{ijk}$ the $B$-quadrea of the triangle with points $%
A_{i}$, $A_{j}$ and $A_{k}$;

\item by $\mathcal{V}$ its $B$-quadrume;

\item by $E_{ij}$ the $B$-dihedral spread between the planes through any
three points intersecting at a line through $A_{i}$ and $A_{j}$;

\item by $\mathcal{S}_{i}$ the $B$-solid spread between three concurrent
lines at $A_{i}$; and

\item by $\mathcal{D}_{i}$ the $B$-dual solid spread between three
concurrent lines at $A_{i}$.
\end{itemize}

The quantity%
\begin{equation*}
\mathcal{R}\equiv \frac{16\mathcal{V}^{2}}{\mathcal{A}_{012}\mathcal{A}_{013}%
\mathcal{A}_{023}\mathcal{A}_{123}}
\end{equation*}%
is a key component of our study, whose geometric meaning is yet to be fully
understood; we will call this the \textbf{Richardson constant}. Here are
some examples of relations we obtain:%
\begin{equation*}
\frac{E_{01}E_{23}}{Q_{01}Q_{23}}=\frac{E_{02}E_{13}}{Q_{02}Q_{13}}=\frac{%
E_{03}E_{12}}{Q_{03}Q_{12}}=\mathcal{R}
\end{equation*}%
\begin{equation*}
\frac{\mathcal{D}_{0}}{\mathcal{A}_{123}}=\frac{\mathcal{D}_{1}}{\mathcal{A}%
_{023}}=\frac{\mathcal{D}_{2}}{\mathcal{A}_{013}}=\frac{\mathcal{D}_{3}}{%
\mathcal{A}_{012}}=\frac{\mathcal{R}}{4}
\end{equation*}%
and%
\begin{eqnarray*}
\frac{\mathcal{S}_{0}\mathcal{S}_{1}\mathcal{S}_{2}}{Q_{03}Q_{13}Q_{23}} &=&%
\frac{\mathcal{S}_{0}\mathcal{S}_{1}\mathcal{S}_{3}}{Q_{02}Q_{12}Q_{23}}=%
\frac{\mathcal{S}_{0}\mathcal{S}_{2}\mathcal{S}_{3}}{Q_{01}Q_{12}Q_{13}}=%
\frac{\mathcal{S}_{1}\mathcal{S}_{2}\mathcal{S}_{3}}{Q_{01}Q_{02}Q_{03}} \\
&=&\frac{\mathcal{V}^{3}}{%
64Q_{01}^{2}Q_{02}^{2}Q_{03}^{2}Q_{12}^{2}Q_{13}^{2}Q_{23}^{2}}.
\end{eqnarray*}%
The first two results are rational/algebraic analogs of some results found
from \cite{Richardson}, while the last one is novel; but we will develop
many more in this paper.

In a later section we specialize our formulas to the case of the
tri-rectangular tetrahedron, and develop some additional key results for
this important situation. We will also use this framework to derive $B$%
-quadrances between opposite lines of a tetrahedron, which are skew, as well
as presenting potential directions for further research.

\section{A review of three-dimensional vector algebra over a general
metrical framework}

We begin with the framework of the three-dimensional vector space $\mathbb{V}%
^{3}$ over a field $\mathbb{F}$ not of characteristic $2$, consisting of row
vectors $v=\left( x,y,z\right) ,$ with the usual arithmetical structures of
vector addition and subtraction, together with scalar multiplication.

\subsection{$B$-scalar product}

A $3\times 3$ symmetric matrix%
\begin{equation}
B\equiv 
\begin{pmatrix}
a_{1} & b_{3} & b_{2} \\ 
b_{3} & a_{2} & b_{1} \\ 
b_{2} & b_{1} & a_{3}%
\end{pmatrix}
\label{SymmetricMatrix}
\end{equation}%
determines a \textbf{symmetric bilinear form}\ on $\mathbb{V}^{3}$ defined by%
\begin{equation*}
v\cdot _{B}w\equiv vBw^{T}\mathbf{.}
\end{equation*}%
We will call this the $B$-\textbf{scalar product}. The associated $B$-%
\textbf{quadratic form }on $\mathbb{V}^{3}$ is defined by\textbf{\ }%
\begin{equation*}
Q_{B}\left( v\right) \equiv v\cdot _{B}v
\end{equation*}%
and we call the number $Q_{B}\left( v\right) $ the $B$\textbf{-quadrance }of%
\textbf{\ }$v$. A vector $v$ is $B$\textbf{-null }precisely when 
\begin{equation*}
Q_{B}\left( v\right) =0.
\end{equation*}%
The $B$-quadrance satisfies the obvious properties that for vectors $v$ and $%
w$ in $\mathbb{V}^{3}$ and a number $\lambda $ in $\mathbb{F}$%
\begin{equation*}
Q_{B}\left( \lambda v\right) =\lambda ^{2}Q_{B}\left( v\right)
\end{equation*}%
as well as%
\begin{equation*}
Q_{B}\left( v+w\right) =Q_{B}\left( v\right) +Q_{B}\left( w\right) +2\left(
v\cdot _{B}w\right)
\end{equation*}%
and 
\begin{equation*}
Q_{B}\left( v-w\right) =Q_{B}\left( v\right) +Q_{B}\left( w\right) -2\left(
v\cdot _{B}w\right) .
\end{equation*}%
Hence the $B$-scalar product can be expressed in terms of the $B$-quadratic
form by either of the two \textbf{polarisation formulas} 
\begin{equation*}
v\cdot _{B}w=\frac{Q_{B}\left( v+w\right) -Q_{B}\left( v\right) -Q_{B}\left(
w\right) }{2}=\frac{Q_{B}\left( v\right) +Q_{B}\left( w\right) -Q_{B}\left(
v-w\right) }{2}.
\end{equation*}%
The $B$-scalar product is \textbf{non-degenerate} precisely when the
condition that $v\cdot _{B}w=0$ for any vector $v$ in $\mathbb{V}^{3}$
implies that $w=0$; this will occur precisely when $B$ is invertible. We
will assume that the $B$-scalar product is non-degenerate throughout this
paper. Finally, two vectors $v$ and $w$ in $\mathbb{V}^{3}$ are $B$-\textbf{%
perpendicular} precisely when $v\cdot _{B}w=0$, in which case we write $%
v\perp _{B}w$.

\subsection{$B$-vector product}

Define the \textbf{adjugate} of $B$ from (\ref{SymmetricMatrix}) to be the
matrix%
\begin{equation*}
\func{adj}B\equiv 
\begin{pmatrix}
a_{2}a_{3}-b_{1}^{2} & b_{1}b_{2}-a_{3}b_{3} & b_{1}b_{3}-a_{2}b_{2} \\ 
b_{1}b_{2}-a_{3}b_{3} & a_{1}a_{3}-b_{2}^{2} & b_{2}b_{3}-a_{1}b_{1} \\ 
b_{1}b_{3}-a_{2}b_{2} & b_{2}b_{3}-a_{1}b_{1} & a_{1}a_{2}-b_{3}^{2}%
\end{pmatrix}%
.
\end{equation*}%
When $B$ is invertible this is 
\begin{equation*}
\func{adj}B=\left( \det B\right) B^{-1}.
\end{equation*}

For vectors $v\equiv \left( v_{1},v_{2},v_{3}\right) $ and $w\equiv \left(
w_{1},w_{2},w_{3}\right) $ in $\mathbb{V}^{3}$, the usual Euclidean vector
product \cite[p. 65]{Gibbs} is 
\begin{equation*}
v\times w\equiv \left(
v_{2}w_{3}-v_{3}w_{2},v_{3}w_{1}-v_{1}w_{3},v_{1}w_{2}-v_{2}w_{1}\right) .
\end{equation*}%
So the $B$\textbf{-vector product} \cite{NotoWild1} of $v$ and $w$ is
defined to be the vector%
\begin{equation*}
v\times _{B}w\equiv \left( v\times w\right) \func{adj}B.
\end{equation*}

With the $B$-scalar and $B$-vector products defined, following \cite%
{NotoWild1} we define the following expressions involving vectors $%
v_{1},v_{2},v_{3}$ and $v_{4}:$

\begin{itemize}
\item $B$\textbf{-scalar triple product}:%
\begin{equation*}
\left[ v_{1},v_{2},v_{3}\right] _{B}\equiv v_{1}\cdot _{B}\left( v_{2}\times
_{B}v_{3}\right)
\end{equation*}

\item $B$\textbf{-vector triple product}:%
\begin{equation*}
\left\langle v_{1},v_{2},v_{3}\right\rangle _{B}\equiv v_{1}\times
_{B}\left( v_{2}\times _{B}v_{3}\right)
\end{equation*}

\item $B$\textbf{-quadruple scalar product}:%
\begin{equation*}
\left[ v_{1},v_{2};v_{3},v_{4}\right] _{B}\equiv \left( v_{1}\times
_{B}v_{2}\right) \cdot _{B}\left( v_{3}\times _{B}v_{4}\right)
\end{equation*}

\item $B$\textbf{-quadruple vector product}:%
\begin{equation*}
\left\langle v_{1},v_{2};v_{3},v_{4}\right\rangle _{B}\equiv \left(
v_{1}\times _{B}v_{2}\right) \times _{B}\left( v_{3}\times _{B}v_{4}\right)
\end{equation*}
\end{itemize}

\subsection{Summary of results of $B$-scalar and vector products}

We summarize some results from\textit{\ }\cite{NotoWild1} pertaining to $B$%
-scalar and $B$-vector products, $B$-triple products and $B$-quadruple
products. The first result allows us to express the $B$-scalar triple
product in terms of determinants.

\begin{theorem}[Scalar triple product theorem]
Let $M\equiv 
\begin{pmatrix}
v_{1} \\ 
v_{2} \\ 
v_{3}%
\end{pmatrix}%
$ for vectors $v_{1}$, $v_{2}$ and $v_{3}$ in $\mathbb{V}^{3}$. Then%
\begin{eqnarray*}
&&\left[ v_{1},v_{2},v_{3}\right] _{B}=\left[ v_{2},v_{3},v_{1}\right] _{B}=%
\left[ v_{3},v_{1},v_{2}\right] _{B} \\
&=&-\left[ v_{1},v_{3},v_{2}\right] _{B}=-\left[ v_{2},v_{1},v_{3}\right]
_{B}=-\left[ v_{3},v_{2},v_{1}\right] _{B} \\
&=&\left( \det B\right) \left[ v_{1},v_{2},v_{3}\right] =\det \left(
MB\right) .
\end{eqnarray*}
\end{theorem}

The following result expresses the $B$-vector triple product as a linear
combination of two vectors.

\begin{theorem}[Lagrange's formula]
For vectors $v_{1}$, $v_{2}$ and $v_{3}$ in $\mathbb{V}^{3}$,%
\begin{equation*}
\left\langle v_{1},v_{2},v_{3}\right\rangle _{B}=\left( \det B\right) \left[
\left( v_{1}\cdot _{B}v_{3}\right) v_{2}-\left( v_{1}\cdot _{B}v_{2}\right)
v_{3}\right] .
\end{equation*}
\end{theorem}

The $B$-scalar quadruple product can be computed as a determinantal identity
involving $B$-scalar products, as follows.

\begin{theorem}[Binet-Cauchy identity]
For vectors $v_{1}$, $v_{2}$, $v_{3}$ and $v_{4}$ in $\mathbb{V}^{3}$,%
\begin{equation*}
\left[ v_{1},v_{2};v_{3},v_{4}\right] _{B}=\left( \det B\right) \left[
\left( v_{1}\cdot _{B}v_{3}\right) \left( v_{2}\cdot _{B}v_{4}\right)
-\left( v_{1}\cdot _{B}v_{4}\right) \left( v_{2}\cdot _{B}v_{3}\right) %
\right] .
\end{equation*}
\end{theorem}

The following result immediately follows from the Binet-Cauchy identity and
links the $B$-vector product between two vectors to their $B$-quadrances and
their $B$-scalar product.

\begin{theorem}[Lagrange's identity]
For vectors $v_{1}$ and $v_{2}$ in $\mathbb{V}^{3}$,%
\begin{equation*}
Q_{B}\left( v_{1}\times _{B}v_{2}\right) =\left( \det B\right) \left[
Q_{B}\left( v_{1}\right) Q_{B}\left( v_{2}\right) -\left( v_{1}\cdot
_{B}v_{2}\right) ^{2}\right] .
\end{equation*}
\end{theorem}

The $B$-vector quadruple product can be computed by using only the $B$%
-scalar triple product, as follows.

\begin{theorem}[Vector quadruple product theorem]
For vectors $v_{1}$, $v_{2}$, $v_{3}$ and $v_{4}$ in $\mathbb{V}^{3}$,%
\begin{eqnarray*}
\left\langle v_{1},v_{2};v_{3},v_{4}\right\rangle _{B} &=&\left( \det
B\right) \left( \left[ v_{1},v_{2},v_{4}\right] _{B}v_{3}-\left[
v_{1},v_{2},v_{3}\right] _{B}v_{4}\right) \\
&=&\left( \det B\right) \left( \left[ v_{1},v_{3},v_{4}\right] _{B}v_{2}-%
\left[ v_{2},v_{3},v_{4}\right] _{B}v_{1}\right) .
\end{eqnarray*}
\end{theorem}

Immediately following from the $B$-vector quadruple product theorem, the
following corollary will prove very useful in the study of a vector
tetrahedron.

\begin{corollary}
Let $M\equiv 
\begin{pmatrix}
v_{1} \\ 
v_{2} \\ 
v_{3}%
\end{pmatrix}%
$ for vectors $v_{1}$, $v_{2}$ and $v_{3}$ in $\mathbb{V}^{3}$. Then%
\begin{equation*}
\left\langle v_{1},v_{2};v_{1},v_{3}\right\rangle _{B}=\left[ \left( \det
B\right) ^{2}\left( \det M\right) \right] v_{1}.
\end{equation*}
\end{corollary}

A consequence of this corollary gives us a result that will also prove
useful in this paper. We present a quick proof of this result as follows.

\begin{theorem}[Scalar triple product of products]
For vectors $v_{1}$,$v_{2}$ and $v_{3}$ in $\mathbb{V}^{3}$, 
\begin{equation*}
\left[ v_{2}\times _{B}v_{3},v_{3}\times _{B}v_{1},v_{1}\times _{B}v_{2}%
\right] _{B}=\left( \det B\right) \left( \left[ v_{1},v_{2},v_{3}\right]
_{B}\right) ^{2}.
\end{equation*}
\end{theorem}

\begin{proof}
For $M\equiv 
\begin{pmatrix}
v_{1} \\ 
v_{2} \\ 
v_{3}%
\end{pmatrix}%
$ we use Corollary 6 and the Scalar triple product theorem to obtain 
\begin{eqnarray*}
\left[ v_{2}\times _{B}v_{3},v_{3}\times _{B}v_{1},v_{1}\times _{B}v_{2}%
\right] _{B} &=&\left[ v_{2}\times _{B}v_{3},v_{1}\times
_{B}v_{2},v_{1}\times _{B}v_{3}\right] _{B} \\
&=&\left[ \left( \det B\right) ^{2}\left( \det M\right) \right] \left[
v_{1},v_{2},v_{3}\right] _{B} \\
&=&\left( \det B\right) \left( \left[ v_{1},v_{2},v_{3}\right] _{B}\right)
^{2}
\end{eqnarray*}%
as required.
\end{proof}

\section{Affine and vector geometry in three dimensions}

For the rest of this paper, we will work over the three-dimensional affine
space over a field $\mathbb{F}$ not of characteristic two, denoted by $%
\mathbb{A}^{3}$, where $\mathbb{V}^{3}$ is its associated vector space.
While the main objects in $\mathbb{A}^{3}$ are \textbf{points}, which we
denote as triples enclosed in rectangular brackets, vectors in $\mathbb{V}%
^{3}$ can be expressed as a separation between two points. In other words,
for two points $X$ and $Y$ in $\mathbb{A}^{3}$, a vector from $X$ to $Y$ is
expressed as $\overrightarrow{XY}$ and computed to be the affine difference $%
Y-X$.

A \textbf{line} is a pair $\left( A,v\right) $ containing a point $A$ in $%
\mathbb{A}^{3}$ and a vector $v$ in $\mathbb{V}^{3},$ so that a point $X$
lies on it precisely when there exists a number $\lambda $ in $\mathbb{F}$
such that%
\begin{equation*}
X-A=\lambda v.
\end{equation*}%
Two lines $\left( A_{1},v_{1}\right) $ and $\left( A_{2},v_{2}\right) $ are
equal precisely when $v_{1}$, $v_{2}$ and $\overrightarrow{A_{1}A_{2}}$ are
all scalar multiples of each other. The vector $v$ is a \textbf{direction
vector }for the line\textbf{\ }$\left( A,v\right) $. Given two points $X_{1}$
and $X_{2}$ both lying on a line, we can denote such a line by $X_{1}X_{2}$,
so that the lines $X_{1}X_{2}$ and $Y_{1}Y_{2}$ are equal precisely when $%
Y_{1}$ and $Y_{2}$ both lie on the line $X_{1}X_{2}$ and vice versa.

A \textbf{plane} in $\mathbb{A}^{3}$\ is a triple $\left( A,v,w\right) $
containing a point $A$ in $\mathbb{A}^{3}$ and two linearly independent
vectors $v$ and $w$ in $\mathbb{V}^{3},$ so that a point $X$ lies on it if
there exists numbers $\lambda $ and $\mu $ in $\mathbb{F}$ such that%
\begin{equation*}
X-A=\lambda v+\mu w.
\end{equation*}%
In other words, the vector $\overrightarrow{AX}$ is a linear combination of $%
v$ and $w$, so that two planes $\left( A_{1},v_{1},w_{1}\right) $ and $%
\left( A_{2},v_{2},w_{2}\right) $ are equal when any of $v_{1}$, $v_{2}$, $%
w_{1}$, $w_{2}$ and $\overrightarrow{A_{1}A_{2}}$ are linear combinations of
any two of these. The vectors $v$ and $w$ are then \textbf{spanning vectors}
for the plane $\left( A,v,w\right) $. We can associate to a plane in $%
\mathbb{A}$ a $B$\textbf{-normal vector} $n$ so that any two points $X$ and $%
Y$ in $\mathbb{A}^{3}$ lying on the plane satisfy%
\begin{equation*}
\left( Y-X\right) \cdot _{B}n=\overrightarrow{XY}\cdot _{B}n=0.
\end{equation*}%
A plane in $\mathbb{A}^{3}$ with three points $X$, $Y$ and $Z$ will be
denoted by $XYZ$, with two planes $X_{1}Y_{1}Z_{1}$ and $X_{2}Y_{2}Z_{2}$
being equal precisely when $X_{2}$, $Y_{2}$ and $Z_{2}$ lie on $%
X_{1}Y_{1}Z_{1}$ and vice versa.

A \textbf{triangle} in $\mathbb{A}^{3}$ is an unordered collection of three
points in $\mathbb{A}^{3}$, say $\left\{ A_{1},A_{2},A_{3}\right\} $, and is
denoted by $\overline{A_{1}A_{2}A_{3}}$. Such a triangle determines two 
\textbf{vector triangles} \cite{NotoWild1} $\left\{ \overrightarrow{%
A_{1}A_{2}},\overrightarrow{A_{2}A_{3}},\overrightarrow{A_{3}A_{1}}\right\} $
and $\left\{ \overrightarrow{A_{1}A_{3}},\overrightarrow{A_{3}A_{2}},%
\overrightarrow{A_{2}A_{1}}\right\} $ where the vectors in the vector
triangle sum to $\mathbf{0}$. By defining $v_{ij}\equiv \overrightarrow{%
A_{i}A_{j}}$ for any integer $i$ and $j$ between $1$ and $3$, we denote
these vector triangles respectively by $\overline{v_{12}v_{23}v_{31}}$ and $%
\overline{v_{13}v_{32}v_{21}}$. A \textbf{tetrahedron} in $\mathbb{A}^{3}$
is an unordered collection of four points in $\mathbb{A}^{3}$, say $\left\{
A_{0},A_{1},A_{2},A_{3}\right\} $, and is denoted by $\overline{%
A_{0}A_{1}A_{2}A_{3}}$. An unordered collection of any two distinct points
of a tetrahedron will be called an \textbf{edge} of the tetrahedron, and an
unordered collection of any three distinct points of a tetrahedron will be
called a \textbf{triangle }of the tetrahedron. Associated to each edge and
triangle of a tetrahedron is the line and plane (respectively) that passes
through the collection of points of the tetrahedron; we call these the 
\textbf{lines} and \textbf{planes} of the tetrahedron.

\section{Affine rational trigonometry in three dimensions}

We now define the rational trigonometric quantities that we will use to
analyse a general tetrahedron over a general field and symmetric bilinear
form.

The $B$-\textbf{quadrance} between two points $A_{1}$ and $A_{2}$ in $%
\mathbb{A}^{3}$ is the number%
\begin{equation*}
Q_{B}\left( A_{1},A_{2}\right) \equiv Q_{B}\left( \overrightarrow{A_{1}A_{2}}%
\right) =\overrightarrow{A_{1}A_{2}}\cdot _{B}\overrightarrow{A_{1}A_{2}}.
\end{equation*}%
Note that $Q_{B}\left( \overrightarrow{A_{1}A_{2}}\right) =Q_{B}\left( 
\overrightarrow{A_{2}A_{1}}\right) $, so that the $B$-quadrance between two
points in $\mathbb{A}^{3}$ is independent of order.

Define\textbf{\ Archimedes' function} \cite[p. 64]{WildDP} as%
\begin{equation*}
A\left( a,b,c\right) \equiv \left( a+b+c\right) ^{2}-2\left(
a^{2}+b^{2}+c^{2}\right)
\end{equation*}%
so that we also have%
\begin{eqnarray*}
A\left( a,b,c\right) &=&4ab-\left( a+b-c\right) ^{2} \\
&=&4ac-\left( a+c-b\right) ^{2} \\
&=&4bc-\left( b+c-a\right) ^{2}.
\end{eqnarray*}%
For a triangle $\overline{A_{1}A_{2}A_{3}}$ with $B$-quadrances%
\begin{equation*}
Q_{1}\equiv Q_{B}\left( A_{2},A_{3}\right) \quad Q_{2}\equiv Q_{B}\left(
A_{1},A_{3}\right) \quad \text{\textrm{and}}\quad Q_{3}\equiv Q_{B}\left(
A_{1},A_{2}\right)
\end{equation*}%
the $B$\textbf{-quadrea} of $\overline{A_{1}A_{2}A_{3}}$ is%
\begin{equation*}
\mathcal{A}_{B}\left( \overline{A_{1}A_{2}A_{3}}\right) \equiv A\left(
Q_{1},Q_{2},Q_{3}\right) .
\end{equation*}%
By the definition of the $B$-quadrance, this is also equal to $\mathcal{A}%
_{B}\left( \overline{v_{12}v_{23}v_{31}}\right) $ and to $\mathcal{A}%
_{B}\left( \overline{v_{13}v_{32}v_{21}}\right) $. So, the $B$-quadrea of a
triangle is simply the $B$-quadrea of either of its two associated vector
triangles.

The following result extends the Quadrea theorem in \cite{NotoWild1} from
the vector triangle setting to the affine triangle setting. As the only
variant to the result is the $B$-quadrea of the affine triangle, we omit the
proof.

\begin{theorem}[Quadrea theorem]
For a triangle $\overline{A_{1}A_{2}A_{3}}$ in $\mathbb{A}^{3}$ with $%
v_{ij}\equiv \overrightarrow{A_{i}A_{j}}$ for any integer $i$ and $j$
between $1$ and $3$, we have 
\begin{eqnarray*}
\frac{\det B}{4}\mathcal{A}_{B}\left( \overline{A_{1}A_{2}A_{3}}\right)
&=&Q_{B}\left( v_{12}\times _{B}v_{31}\right) =Q_{B}\left( v_{12}\times
_{B}v_{23}\right) =Q_{B}\left( v_{23}\times _{B}v_{31}\right) \\
&=&Q_{B}\left( v_{21}\times _{B}v_{13}\right) =Q_{B}\left( v_{21}\times
_{B}v_{32}\right) =Q_{B}\left( v_{32}\times _{B}v_{13}\right) .
\end{eqnarray*}
\end{theorem}

The $B$-\textbf{quadrume} of a tetrahedron $\overline{A_{0}A_{1}A_{2}A_{3}}$
is%
\begin{equation*}
\mathcal{V}_{B}\left( \overline{A_{0}A_{1}A_{2}A_{3}}\right) \equiv \frac{4}{%
\det B}\left[ \overrightarrow{A_{0}A_{1}},\overrightarrow{A_{0}A_{2}},%
\overrightarrow{A_{0}A_{3}}\right] _{B}^{2}.
\end{equation*}%
By the linearity of the $B$-scalar triple product, this will be unchanged if
we base the vectors at another point, for example 
\begin{equation*}
\left[ \overrightarrow{A_{0}A_{1}},\overrightarrow{A_{0}A_{2}},%
\overrightarrow{A_{0}A_{3}}\right] _{B}^{2}=\left[ \overrightarrow{A_{1}A_{0}%
},\overrightarrow{A_{1}A_{2}},\overrightarrow{A_{1}A_{3}}\right] _{B}^{2}.
\end{equation*}%
The following result ensues.

\begin{theorem}[Quadrume product theorem]
If $M\equiv 
\begin{pmatrix}
\overrightarrow{A_{0}A_{1}} \\ 
\overrightarrow{A_{0}A_{2}} \\ 
\overrightarrow{A_{0}A_{3}}%
\end{pmatrix}%
$ then%
\begin{equation*}
\mathcal{V}_{B}\left( \overline{A_{0}A_{1}A_{2}A_{3}}\right) =4\det \left(
MBM^{T}\right) =4\left( \det B\right) \left( \det MM^{T}\right) =4\left(
\det B\right) \left( \det M\right) ^{2}.
\end{equation*}
\end{theorem}

\begin{proof}
The first expression is immediate from the Scalar triple product theorem,
and the others are just rewrites using the multiplicative property of the
determinant.
\end{proof}

The $B$-quadrume is expressed in terms of the $B$-quadrances as follows.

\begin{theorem}[Quadrume theorem]
For a tetrahedron $\overline{A_{0}A_{1}A_{2}A_{3}}$ in $\mathbb{A}^{3}$,
define $Q_{ij}\equiv Q_{B}\left( A_{i},A_{j}\right) $, for integers $i$ and $%
j$ satisfying $0\leq i<j\leq 3$. The $B$-quadrume of the tetrahedron $%
\overline{A_{0}A_{1}A_{2}A_{3}}$ satisfies 
\begin{equation*}
\mathcal{V}_{B}\left( \overline{A_{0}A_{1}A_{2}A_{3}}\right) =\frac{1}{2}%
\begin{vmatrix}
2Q_{01} & Q_{01}+Q_{02}-Q_{12} & Q_{01}+Q_{03}-Q_{13} \\ 
Q_{01}+Q_{02}-Q_{12} & 2Q_{02} & Q_{02}+Q_{03}-Q_{23} \\ 
Q_{01}+Q_{03}-Q_{13} & Q_{02}+Q_{03}-Q_{23} & 2Q_{03}%
\end{vmatrix}%
.
\end{equation*}
\end{theorem}

\begin{proof}
From the Quadrume product theorem,%
\begin{eqnarray*}
\mathcal{V}_{B}\left( \overline{A_{0}A_{1}A_{2}A_{3}}\right) &=&4\det \left(
MBM^{T}\right) \\
&=&4%
\begin{vmatrix}
v_{1}\cdot _{B}v_{1} & v_{1}\cdot _{B}v_{2} & v_{1}\cdot _{B}v_{3} \\ 
v_{1}\cdot _{B}v_{2} & v_{2}\cdot _{B}v_{2} & v_{2}\cdot _{B}v_{3} \\ 
v_{1}\cdot _{B}v_{3} & v_{2}\cdot _{B}v_{3} & v_{3}\cdot _{B}v_{3}%
\end{vmatrix}%
.
\end{eqnarray*}%
By the definition of the $B$-quadratic form and the polarisation formula,
this becomes 
\begin{equation*}
\mathcal{V}_{B}\left( \overline{A_{0}A_{1}A_{2}A_{3}}\right) =\frac{1}{2}%
\begin{vmatrix}
2Q_{01} & Q_{01}+Q_{02}-Q_{12} & Q_{01}+Q_{03}-Q_{13} \\ 
Q_{01}+Q_{02}-Q_{12} & 2Q_{02} & Q_{02}+Q_{03}-Q_{23} \\ 
Q_{01}+Q_{03}-Q_{13} & Q_{02}+Q_{03}-Q_{23} & 2Q_{03}%
\end{vmatrix}%
\end{equation*}
as required.
\end{proof}

The determinant present in the definition is called the \textbf{%
Cayley-Menger determinant} (see \cite{Audet}, \cite[pp. 285-289]{Dorrie} and 
\cite[pp. 124-126]{Sommerville}) and forms a general framework for
calculating higher-dimensional trigonometric quantities of the "distance"
flavour. While named after Cayley and Menger, this formula was known to
Euler and dates back to work of Tartaglia.

Given two lines $l_{1}$ and $l_{2}$ in $\mathbb{A}^{3}$ with respective
direction vectors $v_{1}$ and $v_{2}$, we define the $B$-\textbf{spread}
between them as%
\begin{equation*}
s_{B}\left( l_{1},l_{2}\right) \equiv 1-\frac{\left( v_{1}\cdot
_{B}v_{2}\right) ^{2}}{Q_{B}\left( v_{1}\right) Q_{B}\left( v_{2}\right) }.
\end{equation*}%
Lagrange's identity allows us to rewrite this as%
\begin{equation*}
s_{B}\left( l_{1},l_{2}\right) =\frac{Q_{B}\left( v_{1}\times
_{B}v_{2}\right) }{\left( \det B\right) Q_{B}\left( v_{1}\right) Q_{B}\left(
v_{2}\right) }.
\end{equation*}%
The following result, originally from \cite[p. 82]{WildDP} and proven with $%
B $-vector products in \cite{NotoWild1}, computes the $B$-quadrea of a
triangle in terms of its $B$-quadrances and $B$-spreads. We state it without
proof here.

\begin{theorem}[Quadrea spread theorem]
For a triangle $\overline{A_{1}A_{2}A_{3}}$ with $B$-quadrances%
\begin{equation*}
Q_{1}\equiv Q_{B}\left( A_{2},A_{3}\right) \quad Q_{2}\equiv Q_{B}\left(
A_{1},A_{3}\right) \quad \mathrm{and}\quad Q_{3}\equiv Q_{B}\left(
A_{1},A_{2}\right)
\end{equation*}%
as well as $B$-spreads%
\begin{equation*}
s_{1}\equiv s_{B}\left( A_{1}A_{2},A_{1}A_{3}\right) \quad s_{2}\equiv
s_{B}\left( A_{1}A_{2},A_{2}A_{3}\right) \quad \mathrm{and}\quad s_{3}\equiv
s_{B}\left( A_{1}A_{3},A_{2}A_{3}\right)
\end{equation*}%
and $B$-quadrea $\mathcal{A}\equiv \mathcal{A}_{B}\left( \overline{%
A_{1}A_{2}A_{3}}\right) $, we have that%
\begin{equation*}
\mathcal{A}=4Q_{1}Q_{2}s_{3}=4Q_{1}Q_{3}s_{2}=4Q_{2}Q_{3}s_{1}.
\end{equation*}
\end{theorem}

Given two planes $\Pi _{1}$ and $\Pi _{2}$ in $\mathbb{A}^{3}$ with $B$%
-normal vectors $n_{1}$ and $n_{2}$ respectively, we define the $B$-\textbf{%
dihedral spread} between them to be%
\begin{equation*}
E_{B}\left( \Pi _{1},\Pi _{2}\right) \equiv 1-\frac{\left( n_{1}\cdot
_{B}n_{2}\right) ^{2}}{Q_{B}\left( n_{1}\right) Q_{B}\left( n_{2}\right) }.
\end{equation*}%
This is clearly independent of the rescaling of normal vectors. Note the
similarities between the definition of the $B$-spread and the $B$-dihedral
spread; this is a central theme in projective rational trigonometry (see 
\cite{WildAP} and \cite{WildPS}). The $B$-dihedral spread can also be
rewritten using Lagrange's identity as 
\begin{equation*}
E_{B}\left( \Pi _{1},\Pi _{2}\right) =\frac{Q_{B}\left( n_{1}\times
_{B}n_{2}\right) }{\left( \det B\right) Q_{B}\left( n_{1}\right) Q_{B}\left(
n_{2}\right) }.
\end{equation*}%
The $B$-dihedral spread satisfies the following property.

\begin{theorem}[Dihedral spread theorem]
Let $\Pi _{1}$ be a plane with spanning vectors $v$ and $w_{1}$, and $\Pi
_{2}$ be a plane with spanning vectors $v$ and $w_{2}$, so that these two
planes meet at a line with direction vector $v$. Then,%
\begin{equation*}
E_{B}\left( \Pi _{1},\Pi _{2}\right) =\frac{\left( \det B\right) \left[
v,w_{1},w_{2}\right] _{B}^{2}Q_{B}\left( v\right) }{Q_{B}\left( v\times
_{B}w_{1}\right) Q_{B}\left( v\times _{B}w_{2}\right) }.
\end{equation*}
\end{theorem}

\begin{proof}
We use the rearrangement of the definition of the $B$-dihedral spread using
Lagrange's identity to write%
\begin{equation*}
E_{B}\left( \Pi _{1},\Pi _{2}\right) =\frac{Q_{B}\left( \left( v\times
_{B}w_{1}\right) \times _{B}\left( v\times _{B}w_{2}\right) \right) }{\left(
\det B\right) Q_{B}\left( v\times _{B}w_{1}\right) Q_{B}\left( v\times
_{B}w_{2}\right) }.
\end{equation*}%
By Corollary 6,%
\begin{eqnarray*}
E_{B}\left( \Pi _{1},\Pi _{2}\right) &=&\frac{Q_{B}\left( \left[ \left( \det
B\right) ^{2}\left( \det M\right) \right] v\right) }{\left( \det B\right)
Q_{B}\left( v\times _{B}w_{1}\right) Q_{B}\left( v\times _{B}w_{2}\right) }
\\
&=&\frac{\left( \det B\right) ^{3}\left( \det M\right) ^{2}Q_{B}\left(
v\right) }{Q_{B}\left( v\times _{B}w_{1}\right) Q_{B}\left( v\times
_{B}w_{2}\right) }
\end{eqnarray*}%
where $M\equiv 
\begin{pmatrix}
v \\ 
w_{1} \\ 
w_{2}%
\end{pmatrix}%
$. By the Quadrume product theorem,%
\begin{eqnarray*}
E_{B}\left( \Pi _{1},\Pi _{2}\right) &=&\frac{\left( \det B\right)
^{3}\left( \det M\right) ^{2}Q_{B}\left( v\right) }{Q_{B}\left( v\times
_{B}w_{1}\right) Q_{B}\left( v\times _{B}w_{2}\right) } \\
&=&\frac{\left( \det B\right) \left[ v,w_{1},w_{2}\right] _{B}^{2}Q_{B}%
\left( v\right) }{Q_{B}\left( v\times _{B}w_{1}\right) Q_{B}\left( v\times
_{B}w_{2}\right) }
\end{eqnarray*}%
as required.
\end{proof}

Take three concurrent lines $l_{1}$, $l_{2}$ and $l_{3}$ in $\mathbb{A}^{3}$
with respective direction vectors $v_{1}$, $v_{2}$ and $v_{3}$. We define
the $B$-\textbf{solid spread} between them as%
\begin{equation*}
\mathcal{S}_{B}\left( l_{1},l_{2},l_{3}\right) \equiv \frac{\left( \left[
v_{1},v_{2},v_{3}\right] _{B}\right) ^{2}}{\left( \det B\right) Q_{B}\left(
v_{1}\right) Q_{B}\left( v_{2}\right) Q_{B}\left( v_{3}\right) }.
\end{equation*}%
The $B$-solid spread satisfies the following identity.

\begin{theorem}[Solid spread theorem]
Suppose three lines $l_{1}$, $l_{2}$ and $l_{3}$ in $\mathbb{A}^{3}$ meet at
a single point $O$ with respective direction vectors $v_{1}$, $v_{2}$ and $%
v_{3}$. Furthermore, define the planes%
\begin{equation*}
\Pi _{12}\equiv \left( O,v_{1},v_{2}\right) ,\quad \Pi _{13}\equiv \left(
O,v_{1},v_{3}\right) \quad \mathrm{and\quad }\Pi _{23}\equiv \left(
O,v_{2},v_{3}\right) .
\end{equation*}%
Then,%
\begin{eqnarray*}
\mathcal{S}_{B}\left( l_{1},l_{2},l_{3}\right) &=&E_{B}\left( \Pi _{12},\Pi
_{13}\right) s_{B}\left( l_{1},l_{2}\right) s_{B}\left( l_{1},l_{3}\right) \\
&=&E_{B}\left( \Pi _{12},\Pi _{23}\right) s_{B}\left( l_{1},l_{2}\right)
s_{B}\left( l_{2},l_{3}\right) \\
&=&E_{B}\left( \Pi _{13},\Pi _{23}\right) s_{B}\left( l_{1},l_{3}\right)
s_{B}\left( l_{2},l_{3}\right) .
\end{eqnarray*}
\end{theorem}

\begin{proof}
By rewriting the definition of the $B$-spread using Lagrange's identity, we
have%
\begin{equation*}
s_{B}\left( l_{1},l_{2}\right) =\frac{Q_{B}\left( v_{1}\times
_{B}v_{2}\right) }{\left( \det B\right) Q_{B}\left( v_{1}\right) Q_{B}\left(
v_{2}\right) }\quad \mathrm{and}\quad s_{B}\left( l_{1},l_{3}\right) =\frac{%
Q_{B}\left( v_{1}\times _{B}v_{3}\right) }{\left( \det B\right) Q_{B}\left(
v_{1}\right) Q_{B}\left( v_{3}\right) }.
\end{equation*}%
Given that%
\begin{equation*}
E_{B}\left( \Pi _{12},\Pi _{13}\right) =\frac{\left( \det B\right) \left[
v_{1},v_{2},v_{3}\right] _{B}Q_{B}\left( v_{1}\right) }{Q_{B}\left(
v_{1}\times _{B}v_{2}\right) Q_{B}\left( v_{1}\times _{B}v_{3}\right) }
\end{equation*}%
compute the product of the above three quantities to get our desired result.
The other results follow by symmetry.
\end{proof}

Given three concurrent lines $l_{1}$, $l_{2}$ and $l_{3}$ in $\mathbb{A}^{3}$%
, we construct three lines $k_{12}$, $k_{13}$ and $k_{23}$ with respective
direction vectors%
\begin{equation*}
n_{12}\equiv v_{1}\times _{B}v_{2},\quad n_{13}\equiv v_{1}\times
_{B}v_{3}\quad \mathrm{and}\text{\quad }n_{23}\equiv v_{2}\times _{B}v_{3}
\end{equation*}%
so that all six lines are concurrent and $k_{12}$ is $B$-perpendicular to $%
l_{1}$ and $l_{2}$, $k_{13}$ is $B$-perpendicular to $l_{1}$ and $l_{3},$
and $k_{23}$ is $B$-perpendicular to $l_{2}$ and $l_{3}$. We then define the 
$B$-\textbf{dual solid spread} between lines $l_{1}$, $l_{2}$ and $l_{3}$ to
be%
\begin{equation*}
\mathcal{D}_{B}\left( l_{1},l_{2},l_{3}\right) \equiv \mathcal{S}_{B}\left(
k_{12},k_{13},k_{23}\right) .
\end{equation*}%
We now present an analog to the Solid spread theorem for $B$-dual solid
spreads.

\begin{theorem}[Dual solid spread theorem]
Suppose three lines $l_{1}$, $l_{2}$ and $l_{3}$ in $\mathbb{A}^{3}$ meet at
a single point $O$ with respective direction vectors $v_{1}$, $v_{2}$ and $%
v_{3}$. Furthermore, define the planes%
\begin{equation*}
\Pi _{12}\equiv \left( O,v_{1},v_{2}\right) ,\quad \Pi _{13}\equiv \left(
O,v_{1},v_{3}\right) \quad \mathrm{and\quad }\Pi _{23}\equiv \left(
O,v_{2},v_{3}\right) .
\end{equation*}%
Then,%
\begin{eqnarray*}
\mathcal{D}_{B}\left( l_{1},l_{2},l_{3}\right) &=&s_{B}\left(
l_{1},l_{2}\right) E_{B}\left( \Pi _{12},\Pi _{13}\right) E_{B}\left( \Pi
_{12},\Pi _{23}\right) \\
&=&s_{B}\left( l_{1},l_{3}\right) E_{B}\left( \Pi _{12},\Pi _{13}\right)
E_{B}\left( \Pi _{13},\Pi _{23}\right) \\
&=&s_{B}\left( l_{2},l_{3}\right) E_{B}\left( \Pi _{12},\Pi _{23}\right)
E_{B}\left( \Pi _{13},\Pi _{23}\right) .
\end{eqnarray*}
\end{theorem}

\begin{proof}
First we construct three lines $k_{12}$, $k_{13}$ and $k_{23}$ with
respective direction vectors%
\begin{equation*}
n_{12}\equiv v_{1}\times _{B}v_{2},\quad n_{13}\equiv v_{1}\times
_{B}v_{3}\quad \mathrm{and}\text{\quad }n_{23}\equiv v_{2}\times _{B}v_{3}
\end{equation*}%
so that these three lines are concurrent to $l_{1}$, $l_{2}$ and $l_{3}$. By
the Scalar triple product of products theorem, we know that%
\begin{equation*}
\left[ n_{12},n_{13},n_{23}\right] _{B}=\left( \det B\right) \left( \left[
v_{1},v_{2},v_{3}\right] _{B}\right) ^{2}
\end{equation*}%
so that%
\begin{eqnarray*}
\mathcal{D}_{B}\left( l_{1},l_{2},l_{3}\right) &=&\frac{\left( \left[
n_{12},n_{13},n_{23}\right] _{B}\right) ^{2}}{\left( \det B\right)
Q_{B}\left( n_{12}\right) Q_{B}\left( n_{13}\right) Q_{B}\left(
n_{23}\right) } \\
&=&\frac{\left( \det B\right) \left[ v_{1},v_{2},v_{3}\right] _{B}^{4}}{%
Q_{B}\left( v_{1}\times _{B}v_{2}\right) Q_{B}\left( v_{1}\times
_{B}v_{3}\right) Q_{B}\left( v_{2}\times _{B}v_{3}\right) }.
\end{eqnarray*}%
Now,%
\begin{equation*}
E_{B}\left( \Pi _{12},\Pi _{13}\right) =\frac{\left( \det B\right) \left[
v_{1},v_{2},v_{3}\right] _{B}^{2}Q_{B}\left( v_{1}\right) }{Q_{B}\left(
v_{1}\times _{B}v_{2}\right) Q_{B}\left( v_{1}\times _{B}v_{3}\right) }%
,\quad E_{B}\left( \Pi _{12},\Pi _{23}\right) =\frac{\left( \det B\right) %
\left[ v_{1},v_{2},v_{3}\right] _{B}^{2}Q_{B}\left( v_{2}\right) }{%
Q_{B}\left( v_{1}\times _{B}v_{2}\right) Q_{B}\left( v_{2}\times
_{B}v_{3}\right) }
\end{equation*}%
and, by Lagrange's identity,%
\begin{equation*}
s_{B}\left( l_{1},l_{2}\right) =\frac{Q_{B}\left( v_{1}\times
_{B}v_{2}\right) }{\left( \det B\right) Q_{B}\left( v_{1}\right) Q_{B}\left(
v_{2}\right) }.
\end{equation*}%
Compute the product of the above three quantities to get our desired result.
The other results follow by symmetry.
\end{proof}

\section{Rational trigonometry of a general tetrahedron}

In what follows, we consider a tetrahedron $\overline{A_{0}A_{1}A_{2}A_{3}}$
in $\mathbb{A}^{3}$. For integers $i$ and $j$ satisfying $0\leq i<j\leq 3$,
the $B$-quadrances between any two points $A_{i}$ and $A_{j}$ of $\overline{%
A_{0}A_{1}A_{2}A_{3}}$ will be denoted by $Q_{ij}$; the $B$-quadrea
associated to the triangle $\overline{A_{i}A_{j}A_{k}}$ of $\overline{%
A_{0}A_{1}A_{2}A_{3}}$ will be denoted by $\mathcal{A}_{ijk}$, for integers $%
i$, $j$ and $k$ satisfying $0\leq i<j<k\leq 3$; and its $B$-quadrume will be
denoted by $\mathcal{V}$.

The $B$-spreads between two lines $A_{i}A_{j}$ and $A_{i}A_{k}$ of $%
\overline{A_{0}A_{1}A_{2}A_{3}}$ will be denoted by $s_{i;jk}$, for an
integer $i$ satisfying $0\leq i\leq 3$ and integers $j$ and $k$ distinct
from $i$ satisfying $0\leq j<k\leq 3$, and the $B$-dihedral spreads between
two planes $A_{i}A_{j}A_{k}$ and $A_{i}A_{j}A_{l}$, for integers $i$ and $j$
satisfying $0\leq i<j\leq 3$ and distinct integers $k$ and $l$ between $0$
and $3$ which are also distinct from $i$ and $j$, will be denoted by $E_{ij}$%
.

Finally, the $B$-solid spreads and $B$-dual solid spreads between the three
lines $A_{i}A_{j}$, $A_{i}A_{k}$ and $A_{i}A_{l}$ will be denoted
respectively by $\mathcal{S}_{i}$ and $\mathcal{D}_{i}$, for an integer $i$
satisfying $0\leq i\leq 3,$ and distinct integers $j$, $k$ and $l$ between $%
0 $ and $3$ which are also distinct from $i$.

The Quadrea theorem and Quadrume theorem gives us expressions for $\mathcal{A%
}_{012}$, $\mathcal{A}_{013}$, $\mathcal{A}_{023}$, $\mathcal{A}_{123}$ and $%
\mathcal{V}$ in terms of the six $B$-quadrances of $\overline{%
A_{0}A_{1}A_{2}A_{3}}$. In terms of the $B$-quadrances and $B$-spreads of $%
\overline{A_{0}A_{1}A_{2}A_{3}}$, we use the Quadrea spread theorem to
express the $B$-quadreas as%
\begin{equation*}
\mathcal{A}%
_{012}=4Q_{01}Q_{02}s_{0;12}=4Q_{01}Q_{12}s_{1;02}=4Q_{02}Q_{12}s_{2;01}
\end{equation*}%
\begin{equation*}
\mathcal{A}%
_{013}=4Q_{01}Q_{03}s_{0;13}=4Q_{01}Q_{13}s_{1;03}=4Q_{03}Q_{13}s_{3;01}
\end{equation*}%
\begin{equation*}
\mathcal{A}%
_{023}=4Q_{02}Q_{03}s_{0;23}=4Q_{02}Q_{23}s_{2;03}=4Q_{03}Q_{23}s_{3;02}
\end{equation*}%
and%
\begin{equation*}
\mathcal{A}%
_{123}=4Q_{12}Q_{13}s_{1;23}=4Q_{12}Q_{23}s_{2;13}=4Q_{13}Q_{23}s_{3;12}.
\end{equation*}

\subsection{The Alternating spreads theorem}

The following gives relations between face spreads of a tetrahedron.

\begin{theorem}[Alternating spreads theorem]
For a tetrahedron $\overline{A_{0}A_{1}A_{2}A_{3}}$ with $B$-spreads $%
s_{i;jk}$, for $i$, $j$ and $k$ distinct integers with $0\leq j<k\leq 3$, we
have 
\begin{equation*}
s_{1;02}s_{2;03}s_{3;01}=s_{1;03}s_{2;01}s_{3;02}.
\end{equation*}
\end{theorem}

\begin{proof}
From the Quadrea spread theorem, we know that 
\begin{equation*}
\mathcal{A}_{012}=4Q_{01}Q_{12}s_{1;02}=4Q_{02}Q_{12}s_{2;01}
\end{equation*}%
so that 
\begin{equation*}
\frac{s_{1;02}}{s_{2;01}}=\frac{Q_{02}}{Q_{01}}.
\end{equation*}%
This is also a direct consequence of the Spread law in the triangle $%
\overline{A_{0}A_{1}A_{2}}.$ Similarly we have the relations%
\begin{equation*}
\frac{s_{2;03}}{s_{3;02}}=\frac{Q_{03}}{Q_{02}}\qquad \mathrm{and}\qquad 
\frac{s_{3;01}}{s_{1;03}}=\frac{Q_{01}}{Q_{03}}.
\end{equation*}%
The required result follows by taking the product of these three relations
and cancelling all of the quadrances.
\end{proof}

Note that all the $B$-spreads in the formula involve the index $0$ on the
right hand side; so, including this relation, three other relations hold
which will correspond to the other points of the tetrahedron.

\subsection{Results for $B$-dihedral spreads}

The following result establishes a formula for the $B$-dihedral spread of a
tetrahedron in terms of its $B$-quadrances, $B$-quadreas and $B$-quadrume.

\begin{theorem}[Tetrahedron dihedral spread formula]
For a tetrahedron $\overline{A_{0}A_{1}A_{2}A_{3}}$ with $B$-quadrances $%
Q_{ij}$ for $0\leq i<j\leq 3$, $B$-quadreas $\mathcal{A}_{012}$, $\mathcal{A}%
_{013}$, $\mathcal{A}_{023}$ and $\mathcal{A}_{123}$, and $B$-quadrume $%
\mathcal{V}$, the $B$-dihedral spread $E_{01}$ can be expressed as%
\begin{equation*}
E_{01}=\frac{4Q_{01}\mathcal{V}}{\mathcal{A}_{012}\mathcal{A}_{013}}.
\end{equation*}
\end{theorem}

\begin{proof}
Let%
\begin{equation*}
v_{1}\equiv \overrightarrow{A_{0}A_{1}},\quad v_{2}\equiv \overrightarrow{%
A_{0}A_{2}}\quad \mathrm{and}\quad v_{3}\equiv \overrightarrow{A_{0}A_{3}}.
\end{equation*}%
By the Dihedral spread theorem%
\begin{equation*}
E_{01}=\frac{\left( \det B\right) \left[ v_{1},v_{2},v_{3}\right]
_{B}^{2}Q_{B}\left( v_{1}\right) }{Q_{B}\left( v_{1}\times _{B}v_{2}\right)
Q_{B}\left( v_{1}\times _{B}v_{3}\right) }.
\end{equation*}%
By the Quadrea theorem%
\begin{equation*}
Q_{B}\left( v_{1}\times _{B}v_{2}\right) =\frac{\det B}{4}\mathcal{A}%
_{B}\left( \overline{A_{0}A_{1}A_{2}}\right) =\frac{\det B}{4}\mathcal{A}%
_{012}
\end{equation*}%
and%
\begin{equation*}
Q_{B}\left( v_{1}\times _{B}v_{3}\right) =\frac{\det B}{4}\mathcal{A}%
_{B}\left( \overline{A_{0}A_{1}A_{3}}\right) =\frac{\det B}{4}\mathcal{A}%
_{013}.
\end{equation*}%
Given that 
\begin{equation*}
\mathcal{V}=\frac{4}{\det B}\left[ v_{1},v_{2},v_{3}\right] _{B}^{2}
\end{equation*}%
we combine the above results to get%
\begin{equation*}
E_{01}=\frac{4Q_{01}\mathcal{V}}{\mathcal{A}_{012}\mathcal{A}_{013}}
\end{equation*}%
as required.
\end{proof}

Similarly we have that 
\begin{equation*}
E_{02}=\frac{4Q_{02}\mathcal{V}}{\mathcal{A}_{012}\mathcal{A}_{023}},\quad
E_{03}=\frac{4Q_{03}\mathcal{V}}{\mathcal{A}_{013}\mathcal{A}_{023}},
\end{equation*}%
\begin{equation*}
E_{12}=\frac{4Q_{12}\mathcal{V}}{\mathcal{A}_{012}\mathcal{A}_{123}},\quad
E_{13}=\frac{4Q_{13}\mathcal{V}}{\mathcal{A}_{013}\mathcal{A}_{123}}\quad 
\mathrm{and\quad }E_{23}=\frac{4Q_{23}\mathcal{V}}{\mathcal{A}_{023}\mathcal{%
A}_{123}}.
\end{equation*}%
The following result, a rational version of a result from \cite{Richardson},
allows us to form a relationship between the products of opposite $B$%
-dihedral spreads and the products of opposite $B$-quadrances. For somewhat
mysterious reasons, the quantity%
\begin{equation*}
\mathcal{R}\equiv \frac{16\mathcal{V}^{2}}{\mathcal{A}_{012}\mathcal{A}_{013}%
\mathcal{A}_{023}\mathcal{A}_{123}}
\end{equation*}%
is of significance in the study of the rational trigonometry of a
tetrahedron. Unfortunately we do not currently have a good geometric
interpretation of this quantity, although the two-dimensional analog is the
quadratic curvature of the circumcircle of a triangle, and the following
theorem does provide a partial answer. The quantity $\mathcal{R}$ is the
rational equivalent of a quantity denoted by $h^{2}$ in \cite{Richardson},
and hence we will call it the \textbf{Richardson constant}.

\begin{theorem}[Dihedral spread ratio theorem]
For a tetrahedron $\overline{A_{0}A_{1}A_{2}A_{3}}$ with $B$-quadrances $%
Q_{ij}$ for $0\leq i<j\leq 3$, $B$-quadreas $\mathcal{A}_{012}$, $\mathcal{A}%
_{013}$, $\mathcal{A}_{023}$ and $\mathcal{A}_{123}$, $B$-quadrume $\mathcal{%
V}$ and $B$-dihedral spreads $E_{ij}$ for $0\leq i<j\leq 3$, we have%
\begin{equation*}
\frac{E_{01}E_{23}}{Q_{01}Q_{23}}=\frac{E_{02}E_{13}}{Q_{02}Q_{13}}=\frac{%
E_{03}E_{12}}{Q_{03}Q_{12}}=\mathcal{R}.
\end{equation*}
\end{theorem}

\begin{proof}
From the equivalent formulation of the Dihedral spread theorem for $%
\overline{A_{0}A_{1}A_{2}A_{3}}$, we have%
\begin{equation*}
E_{01}E_{23}=\frac{4Q_{01}\mathcal{V}}{\mathcal{A}_{012}\mathcal{A}_{013}}%
\frac{4Q_{23}\mathcal{V}}{\mathcal{A}_{023}\mathcal{A}_{123}}=\mathcal{R}%
Q_{01}Q_{23}
\end{equation*}%
and similarly%
\begin{equation*}
E_{02}E_{13}=\mathcal{R}Q_{02}Q_{13}\quad \mathrm{and\quad }E_{03}E_{12}=%
\mathcal{R}Q_{03}Q_{12}.
\end{equation*}%
Divide each result through by $Q_{01}Q_{23}$, $Q_{02}Q_{13}$ and $%
Q_{03}Q_{12}$ respectively to obtain our desired result.
\end{proof}

\subsection{Results for $B$-solid spreads}

We now present a formula for calculating the $B$-solid spreads of a
tetrahedron in terms of its $B$-quadrances and $B$-quadrume.

\begin{theorem}[Tetrahedron solid spread formula]
For a tetrahedron $\overline{A_{0}A_{1}A_{2}A_{3}}$ with $B$-quadrances $%
Q_{ij}$ for $0\leq i<j\leq 3$ and $B$-quadrume $\mathcal{V}$, the $B$-solid
spread $\mathcal{S}_{0}$ can be expressed as%
\begin{equation*}
\mathcal{S}_{0}=\frac{\mathcal{V}}{4Q_{01}Q_{02}Q_{03}}.
\end{equation*}
\end{theorem}

\begin{proof}
Defining $v_{1}\equiv \overrightarrow{A_{0}A_{1}},v_{2}\equiv 
\overrightarrow{A_{0}A_{2}}$ and $v_{3}\equiv \overrightarrow{A_{0}A_{3}}$,
the definition of the $B$-solid spread 
\begin{equation*}
\mathcal{S}_{0}=\frac{\left( \left[ v_{1},v_{2},v_{3}\right] _{B}\right) ^{2}%
}{\left( \det B\right) Q_{B}\left( v_{1}\right) Q_{B}\left( v_{2}\right)
Q_{B}\left( v_{3}\right) }.
\end{equation*}%
can be rewritten using the formula for the quadrume $\mathcal{V}$ as%
\begin{equation*}
\mathcal{S}_{0}=\frac{\mathcal{V}}{4Q_{01}Q_{02}Q_{03}}
\end{equation*}%
which is our desired result.
\end{proof}

Similarly, we have%
\begin{equation*}
\mathcal{S}_{1}=\frac{\mathcal{V}}{4Q_{01}Q_{12}Q_{13}},\quad \mathcal{S}%
_{2}=\frac{\mathcal{V}}{4Q_{02}Q_{12}Q_{23}}\quad \mathrm{and\quad }\mathcal{%
S}_{3}=\frac{\mathcal{V}}{4Q_{03}Q_{13}Q_{23}}.
\end{equation*}

We present an interesting result regarding the ratio of $B$-solid spreads.

\begin{theorem}[First solid spread ratio theorem]
For a tetrahedron $\overline{A_{0}A_{1}A_{2}A_{3}}$ with $B$-quadrances $%
Q_{ij}$ for $0\leq i<j\leq 3$, $B$-quadrume $\mathcal{V}$, and $B$-solid
spreads $\mathcal{S}_{k}$ for $0\leq k\leq 3$, we have%
\begin{equation*}
\frac{\mathcal{S}_{0}}{\mathcal{S}_{1}}=\frac{Q_{12}Q_{13}}{Q_{02}Q_{03}}.
\end{equation*}
\end{theorem}

\begin{proof}
This is an immediate consequence of the Tetrahedron solid spread formula, as 
\begin{equation*}
\frac{\mathcal{S}_{0}}{\mathcal{S}_{1}}=\left( \frac{\mathcal{V}}{%
4Q_{01}Q_{02}Q_{03}}\right) \div \left( \frac{\mathcal{V}}{%
4Q_{01}Q_{12}Q_{13}}\right) =\frac{Q_{12}Q_{13}}{Q_{02}Q_{03}}.
\end{equation*}
\end{proof}

Similar formulas hold also for other ratios, for example 
\begin{equation*}
\frac{\mathcal{S}_{0}}{\mathcal{S}_{2}}=\frac{Q_{12}Q_{23}}{Q_{01}Q_{03}}%
\qquad \mathrm{and}\qquad \frac{\mathcal{S}_{1}}{\mathcal{S}_{3}}=\frac{%
Q_{03}Q_{23}}{Q_{01}Q_{12}}\qquad \mathrm{etc}.
\end{equation*}

\begin{theorem}[Second solid spread ratio theorem]
For a tetrahedron $\overline{A_{0}A_{1}A_{2}A_{3}}$ with $B$-quadrances $%
Q_{ij}$ for $0\leq i<j\leq 3$, $B$-quadrume $\mathcal{V}$, and $B$-solid
spreads $\mathcal{S}_{k}$ for $0\leq k\leq 3$ we have%
\begin{equation*}
\frac{\mathcal{S}_{0}\mathcal{S}_{1}}{\mathcal{S}_{2}\mathcal{S}_{3}}=\frac{%
Q_{23}^{2}}{Q_{01}^{2}}.
\end{equation*}
\end{theorem}

\begin{proof}
This is an immediate consequence of the First solid spread ratio theorem, as 
\begin{equation*}
\frac{\mathcal{S}_{0}\mathcal{S}_{1}}{\mathcal{S}_{2}\mathcal{S}_{3}}=\left( 
\frac{\mathcal{S}_{0}}{\mathcal{S}_{2}}\right) \left( \frac{\mathcal{S}_{1}}{%
\mathcal{S}_{3}}\right) =\left( \frac{Q_{12}Q_{23}}{Q_{01}Q_{03}}\right)
\left( \frac{Q_{03}Q_{23}}{Q_{01}Q_{12}}\right) =\frac{Q_{23}^{2}}{Q_{01}^{2}%
}.
\end{equation*}
\end{proof}

Similarly, we will have%
\begin{equation*}
\frac{\mathcal{S}_{0}\mathcal{S}_{2}}{\mathcal{S}_{1}\mathcal{S}_{3}}=\frac{%
Q_{13}^{2}}{Q_{02}^{2}}\quad \mathrm{and}\quad \frac{\mathcal{S}_{0}\mathcal{%
S}_{3}}{\mathcal{S}_{1}\mathcal{S}_{2}}=\frac{Q_{12}^{2}}{Q_{03}^{2}}.
\end{equation*}

We may also derive a result pertaining to the ratio of the product of three $%
B$-solid spreads to the product of three $B$-quadrances; this is a new
result that is unique to this paper, which can only be understood by using
the framework of rational trigonometry.

\begin{theorem}[Third solid spread ratio theorem]
For a tetrahedron $\overline{A_{0}A_{1}A_{2}A_{3}}$ with $B$-quadrances $%
\mathcal{A}_{012}$, $\mathcal{A}_{013}$, $\mathcal{A}_{023}$ and $\mathcal{A}%
_{123}$, $B$-quadrume $\mathcal{V}$, $B$-solid spreads $\mathcal{S}_{0}$, $%
\mathcal{S}_{1}$, $\mathcal{S}_{2}$ and $\mathcal{S}_{3}$, we have%
\begin{equation*}
\frac{\mathcal{S}_{0}\mathcal{S}_{1}\mathcal{S}_{2}}{Q_{03}Q_{13}Q_{23}}=%
\frac{\mathcal{S}_{0}\mathcal{S}_{1}\mathcal{S}_{3}}{Q_{02}Q_{12}Q_{23}}=%
\frac{\mathcal{S}_{0}\mathcal{S}_{2}\mathcal{S}_{3}}{Q_{01}Q_{12}Q_{13}}=%
\frac{\mathcal{S}_{1}\mathcal{S}_{2}\mathcal{S}_{3}}{Q_{01}Q_{02}Q_{03}}=%
\frac{\mathcal{V}^{3}}{%
64Q_{01}^{2}Q_{02}^{2}Q_{03}^{2}Q_{12}^{2}Q_{13}^{2}Q_{23}^{2}}.
\end{equation*}
\end{theorem}

\begin{proof}
By the equivalent formulation of the Solid spread theorem for $\overline{%
A_{0}A_{1}A_{2}A_{3}}$, we have%
\begin{eqnarray*}
\mathcal{S}_{0}\mathcal{S}_{1}\mathcal{S}_{2} &=&\frac{\mathcal{V}}{%
4Q_{01}Q_{02}Q_{03}}\frac{\mathcal{V}}{4Q_{01}Q_{12}Q_{13}}\frac{\mathcal{V}%
}{4Q_{02}Q_{12}Q_{23}} \\
&=&\left( \frac{\mathcal{V}^{3}}{%
64Q_{01}^{2}Q_{02}^{2}Q_{03}^{2}Q_{12}^{2}Q_{13}^{2}Q_{23}^{2}}\right)
Q_{03}Q_{13}Q_{23}.
\end{eqnarray*}%
Divide both sides by $Q_{03}Q_{13}Q_{23}$ to get our desired result. The
other results follow by symmetry.
\end{proof}

\subsection{Results for $B$-dual solid spreads}

We now present a formula for the $B$-dual solid spread of a tetrahedron in
terms of its $B$-quadreas and $B$-quadrume.

\begin{theorem}[Tetrahedron dual solid spread formula]
For a tetrahedron $\overline{A_{0}A_{1}A_{2}A_{3}}$ with $B$-quadreas $%
\mathcal{A}_{012}$, $\mathcal{A}_{013}$, $\mathcal{A}_{023}$ and $\mathcal{A}%
_{123}$, and $B$-quadrume $\mathcal{V}$, the $B$-dual solid spread $\mathcal{%
D}_{0}$ can be expressed as%
\begin{equation*}
\mathcal{D}_{0}=\frac{4\mathcal{V}^{2}}{\mathcal{A}_{012}\mathcal{A}_{013}%
\mathcal{A}_{023}}.
\end{equation*}
\end{theorem}

\begin{proof}
Using the vectors $v_{1}\equiv \overrightarrow{A_{0}A_{1}},v_{2}\equiv 
\overrightarrow{A_{0}A_{2}}$ and $v_{3}\equiv \overrightarrow{A_{0}A_{3}},$
as we saw in the proof of the Dual solid spread theorem, 
\begin{equation*}
\mathcal{D}_{0}=\frac{\left( \det B\right) \left[ v_{1},v_{2},v_{3}\right]
_{B}^{4}}{Q_{B}\left( v_{1}\times _{B}v_{2}\right) Q_{B}\left( v_{1}\times
_{B}v_{3}\right) Q_{B}\left( v_{2}\times _{B}v_{3}\right) }.
\end{equation*}%
But we know that 
\begin{equation*}
\mathcal{V=V}_{B}\left( \overline{A_{0}A_{1}A_{2}A_{3}}\right) \equiv \frac{4%
}{\det B}\left[ v_{1},v_{2},v_{3}\right] _{B}^{2}.
\end{equation*}%
and that 
\begin{eqnarray*}
Q_{B}\left( v_{1}\times _{B}v_{2}\right) &=&\frac{\det B}{4}\mathcal{A}_{012}
\\
Q_{B}\left( v_{1}\times _{B}v_{3}\right) &=&\frac{\det B}{4}\mathcal{A}_{013}
\\
Q_{B}\left( v_{2}\times _{B}v_{3}\right) &=&\frac{\det B}{4}\mathcal{A}_{023}
\end{eqnarray*}%
so that substituting we get 
\begin{eqnarray*}
\mathcal{D}_{0} &=&\frac{\left( \det B\right) \left( \frac{\det B}{4}\right)
^{2}}{\left( \frac{\det B}{4}\right) ^{3}}\frac{\mathcal{V}^{2}}{\mathcal{A}%
_{012}\mathcal{A}_{013}\mathcal{A}_{023}} \\
&=&\frac{4\mathcal{V}^{2}}{\mathcal{A}_{012}\mathcal{A}_{013}\mathcal{A}%
_{023}}.
\end{eqnarray*}
\end{proof}

Similarly we have%
\begin{equation*}
\mathcal{D}_{1}=\frac{4\mathcal{V}^{2}}{\mathcal{A}_{012}\mathcal{A}_{013}%
\mathcal{A}_{123}},\quad \mathcal{D}_{2}=\frac{4\mathcal{V}^{2}}{\mathcal{A}%
_{012}\mathcal{A}_{023}\mathcal{A}_{123}}\quad \mathrm{and\quad }\mathcal{D}%
_{3}=\frac{4\mathcal{V}^{2}}{\mathcal{A}_{013}\mathcal{A}_{023}\mathcal{A}%
_{123}}.
\end{equation*}

The following result outlines the ratio of $B$-dual solid spreads to $B$%
-quadreas, one that is a rational analog of another classical result from 
\cite{Richardson}; it acts similarly to the Sine law in classical
trigonometry, albeit in a very different context. The Richardson constant $%
\mathcal{R}$ will be present here as well.

\begin{theorem}[Dual solid spread and quadrea ratio theorem]
For a tetrahedron $\overline{A_{0}A_{1}A_{2}A_{3}}$ with $B$-quadreas $%
\mathcal{A}_{012}$, $\mathcal{A}_{013}$, $\mathcal{A}_{023}$ and $\mathcal{A}%
_{123}$, $B$-quadrume $\mathcal{V}$, $B$-dual solid spreads $\mathcal{D}_{0}$%
, $\mathcal{D}_{1}$, $\mathcal{D}_{2}$ and $\mathcal{D}_{3}$, and Richardson
constant $\mathcal{R}$, we have%
\begin{equation*}
\frac{\mathcal{D}_{0}}{\mathcal{A}_{123}}=\frac{\mathcal{D}_{1}}{\mathcal{A}%
_{023}}=\frac{\mathcal{D}_{2}}{\mathcal{A}_{013}}=\frac{\mathcal{D}_{3}}{%
\mathcal{A}_{012}}=\frac{\mathcal{R}}{4}.
\end{equation*}
\end{theorem}

\begin{proof}
By the equivalent formulation of the Dual solid spread theorem for $%
\overline{A_{0}A_{1}A_{2}A_{3}}$, we have%
\begin{equation*}
\mathcal{D}_{0}=\frac{4\mathcal{V}^{2}}{\mathcal{A}_{012}\mathcal{A}_{013}%
\mathcal{A}_{023}}.
\end{equation*}%
Divide through by $\mathcal{A}_{123}$ to get%
\begin{equation*}
\frac{\mathcal{D}_{0}}{\mathcal{A}_{123}}=\frac{4\mathcal{V}^{2}}{\mathcal{A}%
_{012}\mathcal{A}_{013}\mathcal{A}_{023}\mathcal{A}_{123}}=\frac{\mathcal{R}%
}{4}.
\end{equation*}%
The other results follow by symmetry.
\end{proof}

\subsection{Tetrahedron skew quadrance formula}

The familiar formula for the projection of one vector onto another (\cite[p.
206]{AR} and \cite[p. 174]{Strang}) holds also for more general bilinear
forms; we define the $B$\textbf{-projection} of a vector $v$ in the
direction of the vector $u$ as the vector%
\begin{equation*}
\left( \func{proj}_{u}v\right) _{B}\equiv \frac{u\cdot _{B}v}{Q_{B}\left(
u\right) }u.
\end{equation*}

For a tetrahedron $\overline{A_{0}A_{1}A_{2}A_{3}}$ with all the above
quantities defined, Hilbert and Cohn-Vossen \cite[pp. 13-17]{Hilbert}
established in the Euclidean case that the pairs of lines $\left(
A_{0}A_{1},A_{2}A_{3}\right) $, $\left( A_{0}A_{2},A_{1}A_{3}\right) $ and $%
\left( A_{0}A_{3},A_{1}A_{2}\right) $ of $\overline{A_{0}A_{1}A_{2}A_{3}}$
are skew, i.e. their meets do not exist. We now define%
\begin{equation*}
n_{01;23}\equiv \overrightarrow{A_{0}A_{1}}\times _{B}\overrightarrow{%
A_{2}A_{3}},\quad n_{02;13}\equiv \overrightarrow{A_{0}A_{2}}\times _{B}%
\overrightarrow{A_{1}A_{3}}\quad \mathrm{and}\quad n_{03;12}\equiv 
\overrightarrow{A_{0}A_{3}}\times _{B}\overrightarrow{A_{1}A_{2}}
\end{equation*}%
so that we define%
\begin{equation*}
R_{01;23}\equiv Q_{B}\left( \func{proj}_{n_{01;23}}\overrightarrow{%
P_{01}P_{23}}\right) ,\quad R_{02;13}\equiv Q_{B}\left( \func{proj}%
_{n_{02;13}}\overrightarrow{P_{02}P_{13}}\right) \quad \mathrm{and}\quad
R_{03;12}\equiv Q_{B}\left( \func{proj}_{n_{03;12}}\overrightarrow{%
P_{03}P_{12}}\right)
\end{equation*}%
to be the \textbf{skew }$B$-\textbf{quadrances} of $\overline{%
A_{0}A_{1}A_{2}A_{3}}$ associated to the respective pairs of opposing lines $%
\left( A_{0}A_{1},A_{2}A_{3}\right) $, $\left( A_{0}A_{2},A_{1}A_{3}\right) $
and $\left( A_{0}A_{3},A_{1}A_{2}\right) $, where $P_{ij}$ is an arbitrary
point on the line $A_{i}A_{j}$ for integers $i$ and $j$ satisfying $0\leq
i<j\leq 3$. This quantity is independent of the selection of the $P_{ij}$'s,
since if the two lines don't meet, then the points on the line will lie on
separate planes which are parallel.\ 

We establish a formula for the skew $B$-quadrances of a tetrahedron based on
its $B$-quadrances and $B$-quadrume. We use \cite{SH} as inspiration to
prove this result in our framework.

\begin{theorem}[Tetrahedron skew quadrance formula]
For a tetrahedron $\overline{A_{0}A_{1}A_{2}A_{3}}$ with $B$-quadrances $%
Q_{ij}$, $B$-quadrume $\mathcal{V}$, and skew $B$-quadrances $R_{01;23}$, $%
R_{02;13}$ and $R_{03;12}$, we have%
\begin{eqnarray*}
R_{01;23} &=&\frac{\mathcal{V}}{4Q_{01}Q_{23}-\left(
Q_{02}+Q_{13}-Q_{03}-Q_{12}\right) ^{2}} \\
R_{02;13} &=&\frac{\mathcal{V}}{4Q_{02}Q_{13}-\left(
Q_{01}+Q_{23}-Q_{03}-Q_{12}\right) ^{2}}
\end{eqnarray*}%
and%
\begin{equation*}
R_{03;12}=\frac{\mathcal{V}}{4Q_{02}Q_{13}-\left(
Q_{01}+Q_{23}-Q_{03}-Q_{12}\right) ^{2}}.
\end{equation*}
\end{theorem}

\begin{proof}
For integers $i$ satisfying $1\leq i\leq 3$, define vectors $v_{i}\equiv 
\overrightarrow{A_{0}A_{i}}$ so that we may define $n_{01;23}\equiv
v_{1}\times _{B}\left( v_{3}-v_{2}\right) $. By the definition of skew $B$%
-quadrances, we have%
\begin{eqnarray*}
R_{01;23} &=&Q_{B}\left( \func{proj}_{n_{01;23}}\overrightarrow{A_{0}A_{2}}%
\right) =Q_{B}\left( \frac{n_{01;23}\cdot _{B}v_{2}}{Q_{B}\left(
n_{01;23}\right) }n_{01;23}\right) \\
&=&\frac{\left[ \left( v_{1}\times _{B}\left( v_{3}-v_{2}\right) \right)
\cdot _{B}v_{2}\right] ^{2}}{Q_{B}\left( v_{1}\times _{B}\left(
v_{3}-v_{2}\right) \right) }=\frac{\left[ v_{1},v_{3}-v_{2},v_{2}\right]
_{B}^{2}}{Q_{B}\left( v_{1}\times _{B}\left( v_{3}-v_{2}\right) \right) }.
\end{eqnarray*}%
Define $M\equiv 
\begin{pmatrix}
v_{1} \\ 
v_{2} \\ 
v_{3}%
\end{pmatrix}%
$ so we may use the bilinearity properties of the $B$-scalar and $B$-vector
products, as well as the Scalar triple product theorem, to rewrite the
numerator as%
\begin{equation*}
\left[ v_{1},v_{3}-v_{2},v_{2}\right] _{B}^{2}=\left( \left[
v_{1},v_{3},v_{2}\right] _{B}-\left[ v_{1},v_{2},v_{2}\right] _{B}\right)
^{2}=\left[ v_{1},v_{2},v_{3}\right] _{B}^{2}=\left( \det M\det B\right)
^{2}.
\end{equation*}%
Furthermore, the denominator becomes%
\begin{eqnarray*}
Q_{B}\left( v_{1}\times _{B}\left( v_{3}-v_{2}\right) \right) &=&\left( \det
B\right) \left[ Q_{B}\left( v_{1}\right) Q_{B}\left( v_{3}-v_{2}\right)
-\left( v_{1}\cdot _{B}\left( v_{3}-v_{2}\right) \right) ^{2}\right] \\
&=&\left( \det B\right) \left( Q_{01}Q_{23}-\left[ \left( v_{1}\cdot
_{B}v_{3}\right) -\left( v_{1}\cdot _{B}v_{2}\right) \right] ^{2}\right)
\end{eqnarray*}%
by Lagrange's identity. Use the polarisation formula to obtain%
\begin{eqnarray*}
Q_{B}\left( v_{1}\times _{B}\left( v_{3}-v_{2}\right) \right) &=&\left( \det
B\right) \left( Q_{01}Q_{23}-\left[ \left( v_{1}\cdot _{B}v_{3}\right)
-\left( v_{1}\cdot _{B}v_{2}\right) \right] ^{2}\right) \\
&=&\left( \det B\right) \left( Q_{01}Q_{23}-\left[ \frac{Q_{01}+Q_{03}-Q_{13}%
}{2}-\frac{Q_{01}+Q_{02}-Q_{12}}{2}\right] ^{2}\right) \\
&=&\left( \det B\right) \left( Q_{01}Q_{23}-\frac{\left(
Q_{02}+Q_{13}-Q_{03}-Q_{12}\right) ^{2}}{4}\right) \\
&=&\frac{\det B}{4}\left( 4Q_{01}Q_{23}-\left(
Q_{02}+Q_{13}-Q_{03}-Q_{12}\right) ^{2}\right) .
\end{eqnarray*}
Combine the results for the numerator and denominator with the Quadrume
product theorem to get%
\begin{eqnarray*}
R_{01;23} &=&\frac{4\left( \det M\det B\right) ^{2}}{\left( \det B\right)
\left( 4Q_{01}Q_{23}-\left( Q_{02}+Q_{13}-Q_{03}-Q_{12}\right) ^{2}\right) }
\\
&=&\frac{4\left( \det B\right) \left( \det M\right) ^{2}}{%
4Q_{01}Q_{23}-\left( Q_{02}+Q_{13}-Q_{03}-Q_{12}\right) ^{2}} \\
&=&\frac{\mathcal{V}}{4Q_{01}Q_{23}-\left(
Q_{02}+Q_{13}-Q_{03}-Q_{12}\right) ^{2}}
\end{eqnarray*}%
as required. The other results follow by symmetry.
\end{proof}

It is curious to note that the denominator of the Tetrahedron skew quadrance
formula is a rational form of Bretschneider's formula \cite{Bretschneider}
for the quadrea of a general quadrangle (a collection of four coplanar
points) in terms of the six quadrances between any two of its points (see 
\cite{Coolidge}, \cite{Dostor} and \cite[pp. 204-205]{Hobson}).

\section{Tri-rectangular tetrahedron}

To finish we apply the framework devised in this paper to study a
particularly fundamental type of tetrahedron, which is the analog of a right
triangle in the three-dimensional setting. Just as many problems in metrical
planar geometry can be resolved into right triangles, and in spherical or
elliptic geometry Napier's rules highlight the importance of right spherical
or elliptic triangles, so the tri-rectangular tetrahedron plays a special
role in three-dimensional geometry.

We set $\overline{A_{0}A_{1}A_{2}A_{3}}$ to be a tetrahedron in $\mathbb{A}%
^{3}$ with all its trigonometric invariants denoted as above. Introducing%
\begin{equation*}
v_{1}\equiv \overrightarrow{A_{0}A_{1}},\quad v_{2}\equiv \overrightarrow{%
A_{0}A_{2}}\quad \mathrm{and}\quad v_{3}\equiv \overrightarrow{A_{0}A_{3}}
\end{equation*}%
we define $\overline{A_{0}A_{1}A_{2}A_{3}}$ to be $B$-\textbf{tri-rectangular%
} \cite[pp. 91-94]{AC}\ at the point $A_{0}$ precisely when $v_{1}$, $v_{2}$
and $v_{3}$ are all mutually $B$-perpendicular, that is when%
\begin{equation*}
v_{1}\cdot _{B}v_{2}=v_{1}\cdot _{B}v_{3}=v_{2}\cdot _{B}v_{3}=0.
\end{equation*}%
While we can also similarly define a $B$-tri-rectangular tetrahedron at
another point of $\overline{A_{0}A_{1}A_{2}A_{3}}$, we may suppose for the
purposes of this study, and without loss of generality, that the tetrahedron 
$\overline{A_{0}A_{1}A_{2}A_{3}}$ is $B$-tri-rectangular at $A_{0}$.

Then by the definition of the $B$-spread we have%
\begin{equation*}
s_{0;12}=s_{0;13}=s_{0;23}=1.
\end{equation*}%
Furthermore, since $v_{1}$, $v_{2}$ and $v_{3}$ are all mutually $B$%
-perpendicular we deduce that%
\begin{equation*}
E_{01}=E_{02}=E_{03}=1.
\end{equation*}%
By the Solid spread projective theorem, we then obtain $\mathcal{S}_{0}=1$.

Because of this, it is natural to parametrize a $B$-tri-rectangular
tetrahedron $\overline{A_{0}A_{1}A_{2}A_{3}}$ by the quadrances%
\begin{equation*}
Q_{01}\equiv K_{1},\quad Q_{02}\equiv K_{2}\quad \mathrm{and}\quad
Q_{03}\equiv K_{3}.
\end{equation*}%
These quantities represent the $B$-quadrances of $\overline{%
A_{0}A_{1}A_{2}A_{3}}$ emanating from the point $A_{0}$. Doing this, we use
Pythagoras' theorem (see \cite{Heath} and \cite{NotoWild1}) to obtain%
\begin{equation*}
Q_{12}=K_{1}+K_{2},\quad Q_{13}=K_{1}+K_{3}\quad \mathrm{and}\quad
Q_{23}=K_{2}+K_{3}.
\end{equation*}%
By the Quadrume theorem, the $B$-quadrume of $\overline{A_{0}A_{1}A_{2}A_{3}}
$ is%
\begin{equation*}
\mathcal{V}=\frac{1}{2}%
\begin{vmatrix}
2K_{1} & 0 & 0 \\ 
0 & 2K_{2} & 0 \\ 
0 & 0 & 2K_{3}%
\end{vmatrix}%
=4K_{1}K_{2}K_{3}.
\end{equation*}%
Then by the Quadrea spread theorem%
\begin{equation*}
\mathcal{A}_{012}=4Q_{01}Q_{02}=4K_{1}K_{2}
\end{equation*}%
\begin{equation*}
\mathcal{A}_{013}=4Q_{01}Q_{03}=4K_{1}K_{3}
\end{equation*}%
and%
\begin{equation*}
\mathcal{A}_{023}=4Q_{02}Q_{03}=4K_{2}K_{3}
\end{equation*}%
are three of its $B$-quadreas. To obtain $\mathcal{A}_{123}$, we rely on the
following generalization of a classical result from \cite{deGua}, which
provides a parallel to Pythagoras' theorem for $B$-quadreas of a $B$%
-tri-rectangular tetrahedron.

\begin{theorem}[de Gua's theorem]
For a $B$-tri-rectangular tetrahedron $\overline{A_{0}A_{1}A_{2}A_{3}}$ at $%
A_{0}$, we have that 
\begin{equation*}
\mathcal{A}_{123}=\mathcal{A}_{012}+\mathcal{A}_{013}+\mathcal{A}_{023}.
\end{equation*}
\end{theorem}

\begin{proof}
By the definition of the $B$-quadrea,%
\begin{eqnarray*}
\mathcal{A}_{123} &=&A\left( Q_{12},Q_{13},Q_{23}\right) \\
&=&\left( 2\left( K_{1}+K_{2}+K_{3}\right) \right) ^{2}-2\left( \left(
K_{1}+K_{2}\right) ^{2}+\left( K_{1}+K_{3}\right) ^{2}+\left(
K_{2}+K_{3}\right) ^{2}\right) \\
&=&4\left(
K_{1}^{2}+K_{2}^{2}+K_{3}^{2}+2K_{1}K_{2}+2K_{1}K_{3}+2K_{2}K_{3}\right)
-4\left(
K_{1}^{2}+K_{2}^{2}+K_{3}^{2}+K_{1}K_{2}+K_{1}K_{3}+K_{2}K_{3}\right) \\
&=&4K_{1}K_{2}+4K_{1}K_{3}+4K_{2}K_{3} \\
&=&\mathcal{A}_{012}+\mathcal{A}_{013}+\mathcal{A}_{023}
\end{eqnarray*}%
as required.
\end{proof}

The result of the Quadrea spread theorem can be rearranged to obtain the
remaining $B$-spreads, which are%
\begin{equation*}
s_{1;02}=\frac{S_{3}K_{1}K_{2}}{K_{1}\left( K_{1}+K_{2}\right) },\quad
s_{1;03}=\frac{K_{1}K_{3}}{K_{1}\left( K_{1}+K_{3}\right) },\quad s_{1;23}=%
\frac{K_{1}K_{2}+K_{1}K_{3}+K_{2}K_{3}}{\left( K_{1}+K_{2}\right) \left(
K_{1}+K_{3}\right) },
\end{equation*}%
\begin{equation*}
s_{2;01}=\frac{K_{1}K_{2}}{K_{2}\left( K_{1}+K_{2}\right) },\quad s_{2;03}=%
\frac{K_{2}K_{3}}{K_{2}\left( K_{2}+K_{3}\right) },\quad s_{2;13}=\frac{%
K_{1}K_{2}+K_{1}K_{3}+K_{2}K_{3}}{\left( K_{1}+K_{2}\right) \left(
K_{2}+K_{3}\right) },
\end{equation*}%
\begin{equation*}
s_{3;01}=\frac{K_{1}K_{3}}{K_{3}\left( K_{1}+K_{3}\right) },\quad s_{3;02}=%
\frac{K_{2}K_{3}}{K_{3}\left( K_{2}+K_{3}\right) }\quad \mathrm{and}\quad
s_{3;12}=\frac{K_{1}K_{2}+K_{1}K_{3}+K_{2}K_{3}}{\left( K_{1}+K_{3}\right)
\left( K_{2}+K_{3}\right) }.
\end{equation*}%
By the Tetrahedron dihedral spread formula, the remaining $B$-dihedral
spreads are%
\begin{equation*}
E_{12}=\frac{K_{3}\left( K_{1}+K_{2}\right) }{%
K_{1}K_{2}+K_{1}K_{3}+K_{2}K_{3}},\quad E_{13}=\frac{K_{2}\left(
K_{1}+K_{3}\right) }{K_{1}K_{2}+K_{1}K_{3}+K_{2}K_{3}}\quad \mathrm{and}%
\quad E_{23}=\frac{K_{1}\left( K_{2}+K_{3}\right) }{%
K_{1}K_{2}+K_{1}K_{3}+K_{2}K_{3}}.
\end{equation*}%
The following elegant relation between $B$-dihedral spreads then becomes
visible.

\begin{theorem}[Tri-rectangular dihedral spread theorem]
For a $B$-tri-rectangular tetrahedron $\overline{A_{0}A_{1}A_{2}A_{3}}$ at $%
A_{0}$, we have that 
\begin{equation*}
E_{12}+E_{13}+E_{23}=2.
\end{equation*}
\end{theorem}

\begin{proof}
Use the above quantities to immediately obtain our result, as follows:%
\begin{equation*}
E_{12}+E_{13}+E_{23}=\frac{K_{3}\left( K_{1}+K_{2}\right) }{%
K_{1}K_{2}+K_{1}K_{3}+K_{2}K_{3}}+\frac{K_{2}\left( K_{1}+K_{3}\right) }{%
K_{1}K_{2}+K_{1}K_{3}+K_{2}K_{3}}+\frac{K_{1}\left( K_{2}+K_{3}\right) }{%
K_{1}K_{2}+K_{1}K_{3}+K_{2}K_{3}}=\allowbreak 2.
\end{equation*}
\end{proof}

By the Tetrahedron solid spread formula, the remaining $B$-solid spreads are%
\begin{equation*}
\mathcal{S}_{1}=\frac{K_{2}K_{3}}{\left( K_{1}+K_{2}\right) \left(
K_{1}+K_{3}\right) },\quad \mathcal{S}_{2}=\frac{K_{1}K_{3}}{\left(
K_{1}+K_{2}\right) \left( K_{2}+K_{3}\right) }\quad \mathrm{and}\quad 
\mathcal{S}_{3}=\frac{K_{1}K_{2}}{\left( K_{1}+K_{3}\right) \left(
K_{2}+K_{3}\right) }.
\end{equation*}%
Recall that another version of Pythagoras' theorem in the plane is that if $%
\overline{A_{0}A_{1}A_{2}}$ has a right angle at $A_{0}$ then 
\begin{equation*}
s_{B}\left( A_{0}A_{1},A_{1}A_{2}\right) +s_{B}\left(
A_{0}A_{2},A_{1}A_{2}\right) =1.
\end{equation*}%
Here is a three-dimensional extension of this involving solid spreads.

\begin{theorem}[Tri-rectangular solid spread theorem]
For a $B$-tri-rectangular tetrahedron $\overline{A_{0}A_{1}A_{2}A_{3}}$ at $%
A_{0}$, we have that 
\begin{equation*}
\left( 1-\mathcal{S}_{1}-\mathcal{S}_{2}-\mathcal{S}_{3}\right) ^{2}=4%
\mathcal{S}_{1}\mathcal{S}_{2}\mathcal{S}_{3}.
\end{equation*}
\end{theorem}

\begin{proof}
With the values of $\mathcal{S}_{1},\mathcal{S}_{2}$ and $\mathcal{S}_{3}$
above, we have 
\begin{eqnarray*}
&&1-\mathcal{S}_{1}+\mathcal{S}_{2}+\mathcal{S}_{3} \\
&=&1-\frac{K_{2}K_{3}}{\left( K_{1}+K_{2}\right) \left( K_{1}+K_{3}\right) }-%
\frac{K_{1}K_{3}}{\left( K_{1}+K_{2}\right) \left( K_{2}+K_{3}\right) }-%
\frac{K_{1}K_{2}}{\left( K_{1}+K_{3}\right) \left( K_{2}+K_{3}\right) } \\
&=&\allowbreak \frac{2K_{1}K_{2}K_{3}}{\left( K_{1}+K_{2}\right) \left(
K_{1}+K_{3}\right) \left( K_{2}+K_{3}\right) }
\end{eqnarray*}%
so that 
\begin{equation*}
\left( 1-\mathcal{S}_{1}-\mathcal{S}_{2}-\mathcal{S}_{3}\right)
^{2}=\allowbreak \frac{4K_{1}^{2}K_{2}^{2}K_{3}^{2}}{\left(
K_{1}+K_{2}\right) ^{2}\left( K_{1}+K_{3}\right) ^{2}\left(
K_{2}+K_{3}\right) ^{2}}=4\mathcal{S}_{1}\mathcal{S}_{2}\mathcal{S}_{3}
\end{equation*}%
as required.
\end{proof}

It appears interesting to ask if this result extends in some fashion to more
general tetrahedra.

As for the dual solid spreads, the value at $A_{0}$ is 
\begin{equation*}
\mathcal{D}_{0}=1
\end{equation*}%
because the normals to the lines meeting there are the lines themselves.
Note that this is consistent with 
\begin{equation*}
\mathcal{D}_{0}=\frac{4\mathcal{V}^{2}}{\mathcal{A}_{012}\mathcal{A}_{013}%
\mathcal{A}_{023}}=\frac{4\left( 4K_{1}K_{2}K_{3}\right) ^{2}}{\left(
4K_{1}K_{2}\right) \left( 4K_{1}K_{3}\right) \left( 4K_{2}K_{3}\right) }%
=\allowbreak 1
\end{equation*}%
which uses the Tetrahedron dual solid spread formula. Using this same
formula, we also get%
\begin{eqnarray*}
\mathcal{D}_{1} &=&\frac{4\mathcal{V}^{2}}{\mathcal{A}_{012}\mathcal{A}_{013}%
\mathcal{A}_{123}}=\frac{4\left( 4K_{1}K_{2}K_{3}\right) ^{2}}{\left(
4K_{1}K_{2}\right) \left( 4K_{1}K_{3}\right) \left(
4K_{1}K_{2}+4K_{1}K_{3}+4K_{2}K_{3}\right) } \\
&=&\frac{K_{2}K_{3}}{K_{1}K_{2}+K_{1}K_{3}+K_{2}K_{3}}
\end{eqnarray*}%
and similarly 
\begin{equation*}
\mathcal{D}_{2}=\frac{K_{1}K_{3}}{K_{1}K_{2}+K_{1}K_{3}+K_{2}K_{3}}\qquad 
\mathrm{and}\qquad \mathcal{D}_{3}=\frac{K_{1}K_{2}}{%
K_{1}K_{2}+K_{1}K_{3}+K_{2}K_{3}}.
\end{equation*}%
The following result pertaining to $\mathcal{D}_{1}$, $\mathcal{D}_{2}$ and $%
\mathcal{D}_{3}$ then follows.

\begin{theorem}[Tri-rectangular dual solid spread]
For a $B$-tri-rectangular tetrahedron $\overline{A_{0}A_{1}A_{2}A_{3}}$ at $%
A_{0}$, we have that 
\begin{equation*}
\mathcal{D}_{1}+\mathcal{D}_{2}+\mathcal{D}_{3}=1.
\end{equation*}
\end{theorem}

\begin{proof}
With the values of $\mathcal{D}_{1}$, $\mathcal{D}_{2}$ and $\mathcal{D}_{3}$
above, we compute that 
\begin{eqnarray*}
&&\mathcal{D}_{1}+\mathcal{D}_{2}+\mathcal{D}_{3} \\
&=&\frac{K_{2}K_{3}}{K_{1}K_{2}+K_{1}K_{3}+K_{2}K_{3}}+\frac{K_{1}K_{3}}{%
K_{1}K_{2}+K_{1}K_{3}+K_{2}K_{3}}+\frac{K_{1}K_{2}}{%
K_{1}K_{2}+K_{1}K_{3}+K_{2}K_{3}} \\
&=&1.
\end{eqnarray*}
\end{proof}

\section{Further directions}

There is clearly a big step in going from the two-dimensional to the
three-dimensional situation in trigonometry. One of the reasons is simply
that the number of objects can increased considerably; instead of just three
points, three lines and a triangle, we have four points, six lines, four
faces, and a tetrahedron, and so the range of metrical notions must also
expand to include the various configurations that are possible when we
combine these in various ways.

So when we contemplate higher dimensional trigonometry, the situation will
become much more involved even when we restrict to the case of the simplex,
and will also require the addition of higher dimensional invariants. We can
expect algebraic relations from the very simple, as in the previous case of
the tri-rectangular tetrahedron, to the enormously complicated and
intricate. We are really only at the beginning of a comprehensive
understanding of the geometry of space.

\end{document}